\documentclass[11pt,leqno]{article}

\usepackage{epsfig}
\usepackage{amssymb}
\usepackage{amscd}
\usepackage[all]{xy}
\usepackage{epsf}

\usepackage{amsmath,amsfonts,amscd,amssymb,theorem}

\newcommand{\eqto}{\mathrel{\stackrel{\sim}{\to}}}

\long\def\comment#1\endcomment{}

\makeatletter
\begingroup
\gdef\th@dotted{\normalfont\itshape
  \def\@begintheorem##1##2{%
        \item[\hskip\labelsep \theorem@headerfont ##1\ ##2.]}%
\def\@opargbegintheorem##1##2##3{%
   \item[\hskip\labelsep \theorem@headerfont ##1\ ##2\ (##3).]}}
\endgroup
\makeatother

\theoremstyle{dotted}

\newtheorem{theorem}{Theorem}[section]
\newtheorem{lemma}[theorem]{Lemma}

\newtheorem{prop}[theorem]{Proposition}
\newtheorem{corr}[theorem]{Corollary}
\newtheorem{klemma}[theorem]{Key-lemma}

\makeatletter
\begingroup
\gdef\th@upshape{\normalfont
  \def\@begintheorem##1##2{%
        \item[\hskip\labelsep \theorem@headerfont ##1\ ##2.]}%
\def\@opargbegintheorem##1##2##3{%
   \item[\hskip\labelsep \theorem@headerfont ##1\ ##2\ (##3).]}}
\endgroup
\makeatother

\theoremstyle{upshape}

\newtheorem{defn}[theorem]{Definition}
\newtheorem{remark}[theorem]{Remark}
\newtheorem{example}[theorem]{Example}

\makeatletter \@addtoreset{equation}{section} \makeatother

\newcommand{\cntrct}                % contraction with a vector field
{\hspace{2pt}\raisebox{1pt}{\text{$\lrcorner$}}\hspace{2pt}}

\newcommand{\proof}[1][Proof.]{\smallskip\noindent{\em #1}}
\def\endproof{\hfill\ensuremath{\square}\par\medskip}

\def\eqref#1{\thetag{\ref{#1}}}

\let\latexref=\ref
\def\ref#1{{\normalfont{\latexref{#1}}}}

\newcommand{\hdot}{{\:\raisebox{3pt}{\text{\circle*{1.5}}}}}
\newcommand{\mb}{\hdot}
\def\dlim_#1{{\displaystyle\lim_{#1}}^\hdot}

\newcommand{\N}{\mathcal{N}}

\newcommand{\Ext}{\operatorname{Ext}}

\newcommand{\Mor}{\mathrm{Mor}}
\newcommand{\Ob}{\mathrm{Ob}}

\newcommand{\Ab}{\mathsf{Ab}}
\newcommand{\Tetra}{\mathsf{Tetra}}

\newcommand{\opp}{\mathrm{opp}}
\newcommand{\BBar}{\mathrm{Bar}}
\newcommand{\Cobar}{\mathrm{Cobar}}
\newcommand{\GS}{\mathrm{GS}}
\newcommand{\Hom}{\mathrm{Hom}}
\newcommand{\Bimod}{\mathsf{Bimod}}
\newcommand{\Bicomod}{\mathsf{Bicomod}}
\newcommand{\Ind}{\mathsf{Ind}}
\newcommand{\Coind}{\mathsf{Coind}}
\newcommand{\EXT}{\underline{\mathcal{E}xt}}
\newcommand{\NR}{\underline{\mathfrak{EXT}}}

\newcommand{\ddeg}{\mathrm{deg}}
\newcommand{\free}{\mathrm{free}}

\newcommand{\sevafigc}[4]{\begin{figure}[h]\centerline{
 \epsfig{file=#1,width=#2,angle=#3}}
\bigskip\caption{#4}\end{figure}}

\def\wtilde#1{\widetilde{#1}\vphantom{#1}}

\title{Hopf algebras, tetramodules, and $n$-fold monoidal categories}
\author{Boris Shoikhet}
\date{}

\begin{document}\maketitle

\begin{abstract}
This paper is an extended version of my talk given in Z\"{u}rich
during the Conference ``Quantization and Geometry'', March 2-6,
2009. The main results are the following.

1. We construct a 2-fold monoidal structure [BFSV] on the category
$\Tetra(A)$ of tetramodules (also known as Hopf bimodules) over an
associative bialgebra $A$. According to an earlier result of
R.Taillefer [Tai1,2], $\Ext^\mb_{\Tetra(A)}(A,A)$ is equal to the
Gerstenhaber-Schack cohomology [GS] of $A$, which governs the
infinitesimal deformations of the bialgebra structure on $A$.

2. Given an $n$-fold monoidal abelian category $\mathcal{C}$ with a common
unit object $A$ and some mild property (*) formulated in the paper,
we consider the graded vector space
$W^\mb=\Ext^\mb_\mathcal{C}(A,A)$. We prove that $W^\mb$ has a
natural $(n+1)$-algebra structure whose product is the Yoneda
product. As a conclusion, the Gerstenhaber-Schack cohomology of any
Hopf algebra $A$ is a 3-algebra.

3. We find an operad of $\mathbb{Z}$-modules which acts on the
Hochschild cohomology of any associative $\mathbb{Z}$-algebra $A$ flat
over $\mathbb{Z}$. The $k$-th component of
this operad is the graded space $\oplus_i\pi^{stab}_{-i}(D^2_k)$ of stable homotopy groups of the
space $D_k^2$, the $k$-th component of the little disks operad. We
establish as well an $n$-monoidal version of this result.

4. We define a contravariant functor from the homotopical category of topological spaces with values in graded vector spaces, depending on an $n$-fold monoidal category (``Hochschild cohomology depending on topological space''). In particular, such a functor is assigned to any associative algebra, any bialgebra, etc.
\end{abstract}
\newpage
\tableofcontents

\section*{Introduction}

\subsection{}
The author would like to warn the reader that this paper is a draft preliminary version; the proofs of some statements in Section 5 are just sketched or even omitted.
\subsection{}
Let $A$ be an associative algebra over a ground field $k$ of characteristic 0, and denote by $HH^\mb(A,A)$ its Hochschild cohomology. It is known that
the graded space $HH^\mb(A,A)$ is a 2-algebra, that is it has a commutative
product of degree 0, a Lie bracket of degree -1, which are
compatible as
\begin{equation}\label{eq1.1}
[a,b\cdot c]=[a,b]\cdot c \pm b\cdot [a,c]
\end{equation}
where the sign $\pm$ in \eqref{eq1.1} is such that the bracket is
{\it odd}, that is for homogeneous $a,b,c$ the sign $\pm=(-1)^{\deg a\deg b +1}$.

These two structures exist also on the Hochschild complex of $A$, where the product is the cup-product of cochains, and the
bracket is the Gerstenhaber bracket. However, on the level of cochains the equation \eqref{eq1.1} fails.

Both the cup-product and the Gerstenhaber bracket are quite
artificial constructions. The cohomology $HH^\mb(A,A)$ is defined as
\begin{equation}\label{eq1.2}
HH^\mb(A,A)=\Ext^\mb_{\Bimod(A)}(A,A)
\end{equation}
where $\Bimod(A)$ is the abelian category of $A$-bimodules, and $A$
is the tautological $A$-bimodule.

It is well-known that the product in \eqref{eq1.1} is the Yoneda product in the definition \eqref{eq1.2}.

Here the following two questions arise: 1) why the Yoneda product is (graded) commutative (for a general abelian category
it is not), and 2) how to define the bracket from the definition \eqref{eq1.2}?

The answer for the both questions uses the fact that $\Bimod(A)$ is
a monoidal category, and $A$ is the unit object in it. The monoidal
structure is the tensor product over $A$, $M_1\otimes_A M_2$. There
is a general theorem which roughly says: given a monoidal category
$\mathcal{C}$ and a unit object $A$ in it,
$\Ext^\mb_\mathcal{C}(A,A)$ is a 2-algebra. (More precisely, there
is a condition (*) on $\mathcal{C}$ formulated later, which is
required for this theorem). This principle was known to many
people, although the author does not know whether any written proof
existed before.

Here we prove this general theorem among other things.

As an example for the general theory, let us prove here that for the unit object $A$ in a monoidal category $\mathcal{C}$ (not necessarily abelian), the
monoid $\Mor_\mathcal{C}(A,A)$ is {\it commutative}.

For denote the monoidal bifunctor by $F\colon\mathcal{C}\times\mathcal{C}\to\mathcal{C}$. Let $A$ be the unit object, and let $f,g\in\Mor(A,A)$.
Then $g\circ f=F\left((id\times g)\circ (f\times id)\right)$ and $f\circ g=F\left((f\times id)\circ (id \times g)\right)$.
But the right-hand sides of the both expressions are equal to $F(f\times g)$, therefore, the left-hand sides are equal.
We proved that $f\circ g=g\circ f\in \Mor(A,A)$.

Our proof that the algebra $\Ext^\mb_\mathcal{C}(A,A)$ is graded-commutative (in the case when $\mathcal{C}$ is a monoidal abelian category) can
be considered as a generalization of the above proof, see Section 2.1.

\subsection{}
Another question on Hochschild cohomology we address here is the following.
Suppose the algebra $A$ is defined over $\mathbb{Z}$, and is flat (= torsion free) over $\mathbb{Z}$.
What is the operad acting in this case on $HH^\mb(A,A)$?

Recall that the operad of 2-algebras can be defined as the homology operad of the little discs operad $\{D^2_n\}$.
That is, for the operad $\mathcal{O}$ governing the 2-algebras, one has:
\begin{equation}\label{intro1}
\mathcal{O}_n=H_\mb(D^2_n;k)
\end{equation}

We prove here the following theorem.
\begin{theorem}
Let $A$ be an algebra over $\mathbb{Z}$ flat over $\mathbb{Z}$. There is a natural action of the operad $\{\mathcal{O}^{\mathbb{Z}}_n\}$
on $HH^\mb(A,A)$, where
\begin{equation}\label{intro2}
\mathcal{O}^\mathbb{Z}_n=\pi_\mb^{stab}(D^2_n)
\end{equation}
\end{theorem}
Here
$\pi^{stab}_k$ stands for the stable homotopy groups,
\begin{equation}\label{intro3}
\pi^{stab}_k(X)=\lim_{s\to\infty}\pi_{k+s}(\Sigma^sX)
\end{equation}
(where $\Sigma$ is the suspension operator, and the limit is attained for any finite CW complex $X$ by the Freudenthal theorem).

The case when $A$ is defined over a field of characteristic 0 is obtained from this more general case of $\mathbb{Z}$-algebra $A$
by the following well-known theorem (see, e.g. [A], Part III, Lecture 4): for any CW complex $X$ one has a canonical (Hurewicz) isomorphism
\begin{equation}\label{intro4}
\pi_k^{stab}(X)\otimes\mathbb{Q}\eqto H_k(X,\mathbb{Q})
\end{equation}
Thus, we recover the 2-algebra operad.

\subsection{}
Another subject of this paper, which originally motivated the overall project, is a proof of following conjecture due to Maxim Kontsevich:
for any associative bialgebra $A$ its Gerstenhaber-Schack cohomology $H^\mb_\GS(A,A)$ is a 3-algebra. The latter means that the graded space $H^\mb_\GS(A,A)$ admits a graded commutative product of degree 0, a Lie bracket of degree -2 such that
\begin{equation}\label{intro5}
[a,b\cdot c]=[a,b]\cdot c \pm b\cdot [a,c]
\end{equation}
where the bracket is {\it even}, that is, for homogeneous $a,b,c$, the sign $\pm=(-1)^{\deg a\deg b}$.

The author tried to use the construction of Stefan Schwede [Sch], which uses the monoidal structure on the category of bimodules over an algebra $A$
in the case of Hochschild cohomology. There is a result of R.Taillefer [Tai1,2] stating that for any bialgebra $A$
\begin{equation}\label{intro6}
H^\mb_\GS(A,A)=\Ext^\mb_{\Tetra(A)}(A,A)
\end{equation}
where $\Tetra(A)$ is the category of {\it tetramodules} over $A$, an abelian category associated to $A$ which is parallel to the category of bimodules
in the case when $A$ is algebra.

We construct on the category of tetramodules {\it two} monoidal structures $\otimes_1$ and $\otimes_2$, which are compatible in some rather
complicated way. These compatibilities altogether give a 2-fold monoidal structure in the sense of [BFSV] on $\Tetra(A)$.

There is a concept of $n$-fold monoidal category, introduced in loc.cit.; it is a main technical tool of this paper. It is a category with $n$ {\it ordered} monoidal structures with some compatibilities. There is an operad of categories acting on any $n$-fold monoidal category; the classifying spaces operad of this operad is homotopy equivalent to the $n$-dimensional little discs operad (this fact is proven in [BFSV]).

We use the last fact to prove the following general theorem.
\begin{theorem}
Let $\mathcal{C}$ be an abelian $n$-fold monoidal category satisfying some mild condition (*), see Section 2.1. Let $A$ be an object which is unit object for all $n$ monoidal structures. Then $\Ext^\mb_{\mathcal{C}}(A,A)$ is naturally an $(n+1)$-algebra whose commutative product is the Yoneda product.
\end{theorem}

When $A$ is a Hopf algebra (that is, a bialgebra with an antipode), the category of tetramodules obeys the condition (*).
As an immediate corollary, one has

\begin{corr}
Let $A$ be a Hopf algebra. Then there is a natural 3-algebra structure on the
Gerstenhaber-Schack cohomology $H^\mb_\GS(A,A)$.
\end{corr}

\subsection{}
The paper is organized as follows:

In Section 1 we recall the theory of Vladimir Retakh [R],[NR] of homotopy groups of the categories of extensions. The Schwede's construction [Sch]
is based on the Retakh's theory, as well as all our generalizations of it. We give a rather detailed exposition with complete proofs, basically because the original note [R] is too coincise, and [NR] works with a great generality of the Waldhausen categories;

In Section 2 we recall the mentioned above construction of Stefan Schwede. This construction gives a Lie bracket on $HH^\mb(A,A)$ in the intrinsic terms of the monoidal abelian category of $A$-bimodules. We discuss here the Schwede's tensor product $\otimes_\tau$, the main ingredient of the construction, and give, following [Sch], an intrinsic construction of the Gerstenhaber bracket on the Hochschild cohomology;

Section 3 introduces $n$-fold monoidal categories. The main new result here is that the category of extensions $\bigsqcup_k\EXT^k_\mathcal{C}(A,A)$
in an $n$-fold monoidal abelian category $\mathcal{C}$ with common unit object $A$, is naturally an $(n+1)$-fold monoidal category;

Section 4 contains our results on the category $\Tetra(A)$ of tetramodules over a bialgebra $A$. We construct a 2-fold monoidal category structure on $\Tetra(A)$. The structure maps $\eta_{MNPQ}$ in this category are {\it not} isomorphisms. To the best of our knowledge, it is the first example of $n$-monoidal categories with the unit object for $n>1$ with this property. In particular, our 2-fold monoidal structure is not defined as a braided category. As well, we discuss the condition (*). It turns out that this condition holds automatically when our bialgebra is a Hopf algebra (that is, admits an antipode). In the case of Hopf algebras we give here a construction of Lie bracket of degree -2 on the Gerstenhaber-Schack cohomology.

In Section 5 we deal with spectra. All main results of the paper, in particular Theorems 0.1 and 0.2 above, are proven here.
The main technical point here is to pass from a ``spectrum" of categories with an action of operad of categories to a spectrum of topological spaces, preserving the action of the corresponding operad. A technical point appears: if we do it naively, the operad action does not admit the ``base points''. The main effort in this Section is directed to introduce somehow the base points, preserving in the same time the homotopical type and the operad action. The basepoints are crucial when we use the smash-product in spectra, which was our initial way to think about the problem. The author tried many ways; finitely, he became successful with a construction imitating the free loop space. We do not achieve basepoints in a proper sense, but we replace them by ``based subcategories''.  Here we are very brief sometimes. We hope to improve the exposition in the sequel version.

Section 6 is served as an Appendix. Here we expose the theory of Rachel Taillefer, claiming that $H^\mb_\GS(A)=\Ext^\mb_{\Tetra(A)}(A,A)$. This result is used throughout in Sections 4 and 5. Our exposition is very closed to the original exposition in [Ta1,2], but the presentation is a bit different.
In particular, it is based on a concept of a $(\mathcal{P},\mathcal{Q})$-pair in an abelian category. As a new result, we give a proof of the well-known folklore statement about the Gerstenhaber-Schack cohomology of the free commutative cocommutative bialgebra $A=S(V)$.

\subsubsection*{Acknowledgements} I am very indebted to Stefan Schwede for many
discussions and explanations. Discussions with Dima Kaledin and
Bernhard Keller were very useful for me. I am thankful to Alexey Gorinov, Alexey Davydov, Pavel Etingof, Ya\"{e}l Fregier, Tom
Leinster, Grisha Merzon, Vladimir Retakh, Dima Tamarkin, Vadik Vologodsky for
discussions and references.
I am especially thankful to Amnon Neeman for many explanations on the work [NR].
I would like to thank the organizers of
the Conference ``Quantization and Geometry'' in Z\"{u}rich Giovanni
Felder and Carlo Rossi for the invitation and for their hospitality.

The work was partially supported by the research grant R1F105L15 of
the University of Luxembourg.

\section{The categories of extensions}

\subsection{}
Let $\mathcal{A}$ be an abelian category, and
let $M,N\in\mathrm{Ob}(\mathcal{A})$, and $k\ge 1$ is integral number.
Consider the following category $\EXT^k(M,N)$.

An object of $\EXT^k(M,N)$ is an exact sequence (an extension)
\begin{equation}\label{eq1.3}
0\rightarrow N\rightarrow F_1\rightarrow\dots\rightarrow
F_k\rightarrow M\rightarrow 0
\end{equation}
and a morphism of two extensions is a map of complexes which is
identity on the ends:
\begin{equation}\label{eq1.4}
\xymatrix{ && F_1\ar[r]\ar[dd]&F_2\ar[r]\ar[dd]&\dots\ar[r]\ar[dd]&F_k\ar[dd]\ar[dr]\\
0\ar[r]& N\ar[ru]\ar[rd]&&&&&M\ar[r]&0\\
&&E_1\ar[r]&E_2\ar[r]&\dots\ar[r]&E_k\ar[ur]}
\end{equation}
Each extension \eqref{eq1.3} defines an element in
$\Ext^k_\mathcal{A}(M,N)$, let us recall the construction (see [M]
for details). For a short exact sequence (extension of length 1)
$0\to A\to B\to C\to 0$ the boundary map in the long exact sequence
defines a map $\delta\colon \Hom(C,C)\to \Ext^1(C,A)$, and the
mentioned assignment is $\delta(id)$. In general case, divide the
long exact sequence into short exact sequences and take the
composition of the maps above.

Thus, we have a map $\varphi\colon \mathrm{Ob}(\EXT^k(M,N))\to
\Ext^k(M,N)$. The natural question is: when two different extensions
define the same elements in $\Ext^k(M,N)$?

The answer goes back to Yoneda and is given e.g. in the MacLane's
book [M]. It is as follows.

\begin{lemma}\label{lemma1.1}
The map $\varphi$ is surjective. Let
$\mathcal{E},\mathcal{F}\in\mathrm{Ob}(\EXT^k(M,N))$ be two
extensions. If there is any morphism $q\colon
\mathcal{E}\to\mathcal{F}$ in $\EXT^k(M,N)$, then
$\varphi(\mathcal{E})=\varphi(\mathcal{F})\in \Ext^k(M,N)$.
Conversably, if $\varphi(\mathcal{E})=\varphi(\mathcal{F})$, the two
extensions $\mathcal{E},\mathcal{F}$ can be connected by a zigzag of
morphisms:
\begin{equation}\label{eq1.7}
\xymatrix{&X_1&&X_3 &&X_\ell\\
\mathcal{E}\ar[ur]&&X_2\ar[ul]\ar[ur]&&\dots\ar[ul]\ar[ur]&&\mathcal{F}\ar[ul]}
\end{equation}
\end{lemma}
The proof can be found in the MacLane's book [M],Ch.3.
\endproof

\subsection{V.Retakh's theorem [R]}
The starting point for the work of Vladimir Retakh [R] is just a
reformulation of the previous Lemma in more contemporary terms.

Consider the nerve ${\mathcal{N}}\EXT^k(M,N)$ (it is a
simplicial set, see [May]), and its geometrical realization
$B\EXT^k(M,N)=|\mathcal{N}\EXT^k(M,N)|$ (it is a topological set, see loc.cit.).
The lemma above is clearly equivalent to the following statement:

\begin{corr}\label{corr1.1}
$\pi_0(\mathcal{N}\EXT^k(M,N))\simeq \Ext^k(M,N)$ as abelian groups.
\end{corr}
Here $\pi_0$ is the 0-th homotopy group, that is, the set of linear
connection components. In our case this set is an abelian group, as
follows. The Baer sum of extensions (see [M], Ch.3) gives a functor
$\mathcal{B}\colon \EXT^k(M,N)\times \EXT^k(M,N)\to\EXT^k(M,N)$
which gives a map of simplicial sets
${\N}\EXT^k(M,N)\times{\N}\EXT^k(M,N)\to{\N}\EXT^k(M,N)$.
This gives the map of classifying spaces
$B\EXT^k(M,N)\times B\EXT^k(M,N)\to B\EXT^k(M,N)$ because the
finite limits commute with the geometric realization functor (see
[May]). Finally, we get a map $\pi_0(B\EXT^k(M,N))\times
\pi_0(B\EXT^k(M,N))\to\pi_0(B\EXT^k(M,N))$. The fact that the Baer
sum becomes the usual sum in $\Ext^k(M,N)$ is proven in [M], Ch.3.
\endproof

Retakh [R] computed the higher homotopy groups of the space
$B\EXT^k(M,N)$. He proves the following theorem:

\begin{theorem}\label{theor1}
\begin{itemize}
\item[1.] For $\ell\le k$,
\begin{equation}\label{eq1.8}
\pi_\ell(\EXT^k(M,N))=\Ext^{k-\ell}(M,N)
\end{equation}
For $\ell>k$, $\pi_\ell(\EXT^k(M,N))=0$;
\item[2.] there are natural homotopy equivalences of topological spaces
\begin{equation}\label{eq1.8.2}
B\EXT^{k-1}(M,N)\to\Omega(B\EXT^k(M,N))
\end{equation}
where $\Omega$ is the loop space functor;
\item[3.] the terms of the spectrum of topological spaces in (2.) are direct products of the Eilenberg-MacLane
spaces.
\end{itemize}
\end{theorem}

\begin{remark}
Note here that there is a paper of Alan Robinson [Rob1], which appeared 5 years before the Retakh's paper, and where an analogous theorem for Tor groups was proven. We would like to mention subsequent papers of Robinson, especially the one on the extraordinary derived category [Rob2]. Probably Robinson was the first who tried to deal with ``modules and algebras in spectra'' in 80's, while the category of symmetric spectra appeared in 90's.
\end{remark}

We overview the main ideas of the proof in Section 1.3.

Let us mention an immediate corollary of the Theorem above.

\begin{corr}\label{characteristic}
Suppose the abelian category $\mathcal{A}$ where we take the
extensions is a $k$-linear abelian category for some field $k$. Then
all homotopy groups $\pi_\ell(\EXT^k(M,N))$, which a priori are
abelian groups, are in fact $k$-vector spaces.
\end{corr}

In particular, if $\mathrm{char}(k)=p$, all elements in
$\pi_\ell(\EXT^k(M,N))$ are $p$-torsion, and if
$\mathrm{char}(k)=0$, all elements in $\pi_\ell(\EXT^k(M,N))$ can
be divided for any integral number.

The formulation of the Theorem above assumes the following lemma:

\begin{lemma}\label{allequivalent}
All connection components in $\EXT^k(A,A)$, $k\ge 0$, are
homotopically equivalent.
\end{lemma}
\proof{} Denote by $\emptyset$ the zero extension with respect to
the Baer sum in $\EXT^k(M,N)$. Suppose $X$ is an extension from
another connected component. Denote the connected components of
$\emptyset$ and $X$ by $\EXT^k(M,N)_\emptyset$ and $\EXT^k(M,N)_X$;
we want to prove that these two categories are homotopy equivalent.
There is a functor $F\colon \EXT^k(M,N)_\emptyset\to \EXT^k(M,N)_X$,
$\alpha\mapsto \alpha+X$ (here $+$ is the Baer sum). Consider
$Y\in\EXT^k(M,N)$ such that $X+Y$ is zero element on $\Ext^k(M,N)$.
Then we have a functor $G\colon \EXT^k(M,N)_X\to
\EXT^k(M,N)_\emptyset$, $\beta\mapsto \beta+Y$. There is a zigzag of
morphisms from $X+Y$ to $\emptyset$ in $\EXT^k(M,N)$, let it be
\begin{equation}\label{eq2.7}
\xymatrix{ &C_1\ar[dl]\ar[dr]&&C_3\ar[dl]\ar[dr]&&C_{2s-1}\ar[dl]\ar[dr]\\
\emptyset&&C_2&&\dots&&X+Y}
\end{equation}
All these arrows can be considered as natural transformation of
functors, say from the functor $F_t\colon
\EXT^k(M,N)_\emptyset\to\EXT^k(M,N)_\emptyset$, $\alpha\mapsto
\alpha+C_t$, to the functor $F_{t+1}\colon
\EXT^k(M,N)_\emptyset\to\EXT^k(M,N)_\emptyset$, $\alpha\mapsto
\alpha+ C_{t+1}$. Therefore, the maps of the classifying space of
$\EXT^k(M,N)_\emptyset$ to itself, induced by these functors, are
homotopic. Performing this for all arrows of the zigzag, we get that
$G\circ F$ is homotopic to the identity. Analogously one proves that
$F\circ G$ is homotopic to the identity.
\endproof

\comment

\subsection{A proof of Theorem \ref{theor1}}\label{sec1.4}
Here we give a proof of the Retakh's theorem not just for completeness; it will be used throughout and generalized in the last Section of
the paper.

We follow the paper [NR], where the proof is slightly different from the Retakh's paper [R].
The author is thankful to Amnon Neeman for explanations on [NR] in the simplest case of abelian category.
See [Q], Section 1, for more details on the homotopy theory of categories.
\subsubsection{The plan}
The idea of the proof is to construct a category $\NR^k(M,N)$ with a functor $F\colon \NR^k(M,N)\to\EXT^k(M,N)$
with the following two properties:
\begin{itemize}
\item[(i)] the category $\NR^k(M,N)$ is contractible;
\item[(ii)] the functor $F$ obeys the conditions of the Corollary Quillen's Theorem B;
\item[(iii)] the fiber of the functor $F^{-1}(*)$ over some distinguished object $*\in\EXT^k(M,N)$ is equivalent to
the category $\EXT^{k-1}(M,N)$.
\end{itemize}

It follows from (ii) that the fiber $F^{-1}(*)$ is equivalent to the homotopical fiber; it follows from (i) that classifying space
of the homotopical fiber is the loop space of the classifying space of the base; it follows from (iii) that the latter space is
$\Ext^{k-1}(M,N)\simeq\Omega\mathcal{N}\EXT^{k}(M,N)$.

\subsubsection{The category $\NR^k(M,N)$}
Define the category $\NR^k(M,N)$ as follows. An object of it is the following commutative diagram
\begin{equation}\label{nr}
\xymatrix{
0\ar[r]&N\ar[d]_{id}\ar[r]&X_1\ar[d]\ar[r]&\dots\ar[r]&X_k\ar[d]\ar[r]&M\ar[r]\ar[d]&0\\
0\ar[r]&N\ar[r]&Y_1\ar[r]&\dots\ar[r]&Y_k\ar[r]&0
}
\end{equation}
where the two rows are exact, and the left vertical arrow is identity. A morphism is a morphism of diagrams.

We prove
\begin{lemma}
The category $\NR^k(M,N)$ is contractible.
\end{lemma}
\proof{}
According to [Q], Proposition 2 on page 84, a natural transformation between two functors $F_1,F_2\colon \mathcal{C}_1\to
\mathcal{C}_2$ gives a homotopy between the corresponding maps $F_{1*},F_{2*}\colon\mathcal{N}\mathcal{C}_1\to\mathcal{NC}_2$.

There is a natural transformation from the identity functor $id\colon \NR^k(M,N)\to \NR^k(M,N)$ to
the functor $\Phi\colon \NR^k(M,N)\to\NR^k(M,N)$ which maps the diagram \eqref{nr}
to the diagram
\begin{equation}\label{nrphi}
\xymatrix{0\ar[r]&N\ar[d]_{id}\ar[r]& Y_1\ar[d]_{id}\ar[r]&\dots\ar[r]&Y_{k-1}\ar[d]_{id}\ar[r]&Y_k\oplus M\ar[d]_{id\oplus pr}\ar[r]^{pr}&M\ar[r]\ar[d]&0\\
0\ar[r]&N\ar[r]&Y_1\ar[r]&\dots\ar[r]&Y_{k-1}\ar[r]& Y_k\ar[r]&0
}
\end{equation}

There is another natural transformation to $\Phi$ from a functor $\Psi\colon\NR^k(M,N)\to\NR^k(M,N)$ where the functor $\Psi$ maps the diagram
\eqref{nr} to
\begin{equation}\label{nrpsi}
\xymatrix{0\ar[r]&N\ar[r]^{id}\ar[d]_{id}&N\ar[r]\ar[d]_{id}&0\ar[d]\ar[r]&\dots\ar[r]&0\ar[d]\ar[r]&M\ar[r]^{id}\ar[d]&M\ar[r]\ar[d]&0\\
0\ar[r]&N\ar[r]^{id}&N\ar[r]&0\ar[r]&\dots\ar[r]&0\ar[r]&0\ar[r]&0}
\end{equation}
The two homotopies together contract the identity map of $\mathcal{N}\NR^k(M,N)$ to a point.
\endproof

\subsubsection{The Quillen's Theorem B}
There is a natural projection functor $F\colon \NR^k(M,N)\to\EXT^k(M,N)$. We are going to prove that this functor
satisfies the conditions of the Quillen's Theorem B. Let us recall what it is. Actually we formulate the Corollary to Theorem B which we
will only need.

Let $g\colon E\to B$ be a map of topological spaces, $b\in B$.
The {\it homotopical fiber} $F(g,b)$ of the map $g$ over $b$ is the spaces
\begin{equation}\label{homfib}
F(g,b)=\{e\in E \text{ and a path } p\colon [0,1]\to B|\ p(0)=g(e),\ p(1)=b\}
\end{equation}
In other words, we replace a map $g\colon E\to B$ by a homotopical to it Serre fibration, and then compute the fiber.
If the map $g\colon E\to B$ is a Serre fibration itself, then up to homotopy the homotopical fiber coincides with the set-theoretical
fiber.

\begin{example}
Let $G$ be a topological group considered as a category $\underline{G}$ with a single object, and let $BG$ be its classifying space.
There is a bundle $p\colon EG\to BG$ which is a Serre fibration with fiber $G$, and with $EG$ homotopically trivial. Then $G$ is homotopically
equivalent to the homotopical fiber, which is $\Omega BG$ since $EG$ is contractible.
$$
G\simeq\Omega BG
$$
\end{example}

Let $f\colon \mathcal{C}_1\to\mathcal{C}_2$ be a functor. For any $Y\in\mathrm{Ob}(\mathcal{C}_2)$ define the
{\it comma category} $Y\setminus f$ whose objects are the pairs $\{(X,v)|X\in \mathrm{Ob}\mathcal{C}_1,\ v\in\Mor_{\mathcal{C}_2}(Y,f(X))\}$.
The morphisms are the maps $X\to X^\prime$ such that the corresponding triangle is commutative.

Morally, the comma category $Y\setminus f$ has the role of the fiber category
$f^{-1}(Y)$ consisting from objects over $Y$ and their morphisms which are mapped by $f$ to the identity of $Y$.
Let us discuss a condition under which the categories $Y\setminus f$ and $f^{-1}(Y)$ are indeed homotopically equivalent.

One says that a functor $f\colon \mathcal{C}_1\to\mathcal{C}_2$ makes the category
$\mathcal{C}_1$ a {\it prefibred category} over $\mathcal{C}_2$ if for every object $Y\in\mathrm{Ob}\mathcal{C}_2$ the functor
\begin{equation}\label{prefibred}
f^{-1}(Y)\to Y\setminus f,\ \ \ X\mapsto (X,id_Y)
\end{equation}
has a right adjoint.

In this case it is clear (see [Q], Corollary 1 on page 84) that the functor above is a homotopical equivalence.

In the case above denote by $(X,v)\mapsto v^*X$ the right adjoint functor to\eqref{prefibred}. Then for any arrow $u\colon Y\to Y^\prime$
in $\mathcal{C}_2$ one has the following functor $u^*\colon f^{-1}(Y^\prime)\to f^{-1}(Y)$:
\begin{equation}\label{prefibred2}
f^{-1}(Y^\prime)\rightarrow Y^\prime\setminus f \xrightarrow{u} Y\setminus f \xrightarrow{v^*} f^{-1}(Y)
\end{equation}
This functor is called {\it the base change functor}.

\newpage

Now we are ready to formulate the {\bf Corollary to the Quillen's Theorem B}:
\begin{theorem}\label{quillenB}
Suppose a functor $f\colon \mathcal{C}_1\to\mathcal{C}_2$ is prefibred and that for any arrow $u\colon Y\to Y^\prime$ in $\mathcal{C}_2$
the base change functor $u^\prime\colon f^{-1}(Y^\prime)\to f^{-1}(Y)$ is a homotopy equivalence.
Then for any $Y\in\mathrm{Ob}\mathcal{C}_2$, the classifying space of the category $f^{-1}(Y)$ is homotopy equivalent to the homotopy-fibre of
$f\colon \mathcal{NC}_1\to \mathcal{NC}_2$ over $Y$.
\end{theorem}
See [Q], Section 1 for a proof.
\endproof

\subsubsection{The functor $F\colon\NR^k(M,N)\to\EXT^k(M,N)$ satisfies the hypothesis of Theorem \ref{quillenB}}
Fix $t\in\mathrm{Ob}(\EXT^k(M,N))$ and consider the two categories $F^{-1}(t)$ and $t\setminus F$ (see previous Subsection).
There is a natural inclusion functor $j\colon F^{-1}(t)\to t\setminus F$. We prove
\begin{lemma}
The functor $j$ has a right adjoint $v^*$.
\end{lemma}
\proof{}
Let $t$ be an exact sequence
\begin{equation}\label{t}
0\rightarrow N\rightarrow T_1\rightarrow\dots\rightarrow T_k\rightarrow M\rightarrow 0
\end{equation}
and let $s\in t\setminus F$ be a pair $(S,v)$ where
$S$ is
\begin{equation}\label{setminus}
\xymatrix{
0\ar[r]&N\ar[d]_{id}\ar[r]&X_1\ar[d]\ar[r]&\dots\ar[r]&X_k\ar[d]\ar[r]&M\ar[r]\ar[d]&0\\
0\ar[r]&N\ar[r]&Y_1\ar[r]&\dots\ar[r]&Y_k\ar[r]&0
}
\end{equation}
and $v$ is a morphism of extensions
\begin{equation}\label{v}
\xymatrix{
0\ar[r]&N\ar[d]_{id}\ar[r]&T_1\ar[d]\ar[r]&\dots\ar[r]&T_k\ar[d]\ar[r]&M\ar[d]_{id}\ar[r]&0\\
0\ar[r]&N\ar[r]&X_1\ar[r]&\dots\ar[r]&X_k\ar[r]&M\ar[r]&0}
\end{equation}
Define $v^*S$ as the diagram
\begin{equation}\label{vstars}
\xymatrix{
0\ar[r]&N\ar[d]_{id}\ar[r]&T_1\ar[d]\ar[r]&\dots\ar[r]&T_k\ar[d]\ar[r]&M\ar[r]\ar[d]&0\\
0\ar[r]&N\ar[r]&Y_1\ar[r]&\dots\ar[r]&Y_k\ar[r]&0}
\end{equation}
where the vertical arrows are the compositions of those in \eqref{v} and in \eqref{setminus}.

We claim that
\begin{equation}\label{radjoint}
\Hom_{F^{-1}(t)}(q, v^*S)=\Hom_{t\setminus F}(j(q),S)
\end{equation}
as bifunctors. But this is clear.
\endproof

It follows from this Lemma that the functor $j\colon F^{-1}(t)\to t\setminus F$ is a homotopy equivalence (see previous Subsection).

Now let $u\colon t^\prime \to t$ be a morphism in $\EXT^k(M,N)$. We need to check that the base change functor $u^*\colon
F^{-1}(t)\to F^{-1}(t^\prime)$ is a
homotopy equivalence.

\begin{lemma}\label{base}
For any morphism $u\colon t^\prime \to t$ in $\EXT^k(M,N)$ the base change functor $u^*\colon
F^{-1}(t)\to F^{-1}(t^\prime)$ is a homotopy equivalence.
\end{lemma}
\proof{}
The functor $u^*$ is just the composition of \eqref{vstars} with $u$. We try to construct a right adjoint functor $R$ for it.
We need
\begin{equation}\label{ur}
\Hom_{F^{-1}(t^\prime)}(u^*\alpha,\beta)=\Hom_{F^{-1}(t)}(\alpha, R\beta)
\end{equation}
For $\alpha\in F^{-1}(t)$, $u^*\alpha$ is just the composition with $u$. We have $u\colon T^\prime\to T$ and $\beta\colon T^\prime\to Y^\prime$.
Consider the cocartesian coproduct $T\sqcup_{T^\prime} Y^\prime$, and the corresponding map $\phi\colon T\to T\sqcup_{T^\prime} Y^\prime$.
Define $R\beta$ as this map $\phi$, the diagram \eqref{altproof}.
\begin{equation}\label{altproof}
\xymatrix{
T^\prime\ar[d]_{u}\ar[r]^{\beta}&Y^\prime\ar[d]\\
T\ar[r]^{\phi}\ar[d]_{\alpha}&T\sqcup_{T^\prime}Y^\prime\ar[dl]\\
Y}
\end{equation}
In general, the sequence $T\sqcup_{T^\prime} Y^\prime$ is not exact; it is if the abelian category satisfies the axiom AB$5$ from [Tohoku]
(all filtred coproducts are exact).
All abelian categories we consider in this paper satisfy AB5.
To treat the general case, when AB5 fails, the proof follows when $u\colon T^\prime\to T$ is a degree-wise {\it monomorphism}.
It is possible to reduce the general case to the one, see [NR], Section 5 for detail.
\endproof

\subsubsection{The fiber $F^{-1}(t)$ for some $t\in\EXT^k(M,N)$ is homotopy equivalent to $\EXT^{k-1}(M,N)$}
Suppose now that $k\ge 2$, and consider the following element $t\in \EXT^k(M,N)$:
\begin{equation}\label{zvezda}
0\rightarrow N\xrightarrow{id}N\rightarrow 0\rightarrow 0\rightarrow\dots\rightarrow 0\rightarrow M\xrightarrow{id}M\rightarrow 0
\end{equation}
(the case $k=1$ should be considered separately, see [NR]).

We are going to compute the category $F^{-1}(t)$.
\begin{lemma}
For $k\ge 2$ and $t$ given in \eqref{zvezda}, the category $F^{-1}(t)$ is homotopy equivalent to $\EXT^{k-1}(M,N)$.
\end{lemma}
\proof{}
The category $F^{-1}(t)$ is the category $\mathcal{C}$ of the diagrams
\begin{equation}\label{reduction}
\xymatrix{&&&&&M\ar[d]^{\phi}\\
0\ar[r]&N\ar[r]&Y_1\ar[r]&\dots\ar[r]&Y_{k-1}\ar[r]&Y_k\ar[r]&0}
\end{equation}
where $\phi$ is an arbitrary map, and morphisms are morphisms of diagrams.

The pullback by $\phi$ gives a functor from this category to $\EXT^{k-1}(M,N)$. The identity functor on $\mathcal{C}$ is homotopy equivalent to the pullback functor.
\endproof

\subsubsection{The Retakh's theorem}
After the Lemmas above, the first two claims of the Retakh's Theorem \ref{theor1} easily follow from the Quillen's Theorem B (Theorem \ref{quillenB} above).

See [NR], Section 8, on the third claim.

\endcomment

\subsection{A proof of the Retakh's theorem}
Here we recall the proof [R] of Theorem \ref{theor1}. Some constructions in this proof will be used later in Section 5.
\subsubsection{The idea}
Let $*$ be some fixed object in $\EXT^k(M,N)$. Consider the category of 1-pathes $\mathcal{P}_R\EXT^k(M,N)$. An object of this category is a ``path''
\begin{equation}\label{bis1}
*\rightarrow X\leftarrow Y
\end{equation}
where $X,Y\in\EXT^k(M,N)$, and the arrows are morphisms of extensions. A morphism in the category $\mathcal{P}_R\EXT^k(M,N)$ is defined as a commutative diagram
\begin{equation}\label{bis2}
\xymatrix{
{*}\ar[r]\ar[d]_{id}&X\ar[d]&Y\ar[l]\ar[d]\\
{*}\ar[r]&X^{\prime}&Y^{\prime}\ar[l]}
\end{equation}
There is a natural projection $p\colon \mathcal{P}_R\EXT^k(M,N)\to\EXT^k(M,N)$, which assigns $Y$ to the object \eqref{bis1}.

We prove the following three statements:
\begin{itemize}
\item[1.] the category $\mathcal{P}_R\EXT^k(M,N)$ is contractible;
\item[2.] the projection $p\colon \mathcal{P}_R\EXT^k(M,N)\to\EXT^k(M,N)$ satisfies the hypothesis of the Quillen's Theorem B;
\item[3.] for some explicit choice of the based object $*$, the fiber $p^{-1}(*)$ is homotopy equivalent to the category $\EXT^{k-1}(M,N)$.
\end{itemize}
The second statement means that one can think on the projection $p$ as on a fibration of topological spaces, by the first statement the total space of this fibration is contractible. These two statements together imply that the set-theoretical fiber $p^{-1}(*)$ is homotopy equivalent to the homotopy fiber, which in the case of contractible total space is equal to the loop spaces on the base. Then, together with the third statement, we get that
$B\EXT^{k-1}(M,N)=\Omega(B\EXT^k(M,N))$. This speculation proves the first two claims of Theorem \ref{theor1}. We refer the reader to [NR], Section 8 for a proof of the third claim of Theorem \ref{theor1}.

\subsubsection{The category $\mathcal{P}_R\EXT^k(M,N)$ is contractible}
First of all, let us recall some general principles on the homotopy theory of categories (see [Q], Section 1).

Any functor $F\colon\mathcal{C}_1\to\mathcal{C}_2$ defines a map $F_B\colon B\mathcal{C}_1\to B\mathcal{C}_2$.
\begin{lemma}\label{elem}
\begin{itemize}
\item[(i)] Any natural transformation between two functors $F,G\colon \mathcal{C}_1\to\mathcal{C}_2$ defines a homotopy between the maps
$F_B,G_B\colon B\mathcal{C}_1\to B\mathcal{C}_2$;
\item[(ii)] if a functor $F\colon \mathcal{C}_1\to\mathcal{C}_2$ admits left or right adjoint, the map $F_B$ is a homotopy equivalence;
\item[(iii)] in particular, if a category $\mathcal{C}$ has initial or final object, the space $B\mathcal{C}$ is contractible.
\end{itemize}
\end{lemma}
All these statements are fairly simple, see [Q], Section 1.
\endproof

We say that a category $\mathcal{C}$ is contractible if the topological space $B\mathcal{C}$ is contractible.

We prove
\begin{lemma}
For any (small) category $\mathcal{C}$, the category $\mathcal{P}_R\mathcal{C}$ is contractible.
\end{lemma}
\proof{}
Define the category $\mathcal{C}(*)$ as the category of objects of $\mathcal{C}$ under $*$. That is, an object of $\mathcal{C}(*)$
is a pair $(X,\varphi)$ where $X$ is an object of $\mathcal{C}$ and $\varphi: *\to X$ is a morphism. A morphism in $\mathcal{C}(*)$ is a commutative diagram
\begin{equation}\label{bis5}
\xymatrix{
&X\ar[dd]\\
{*}\ar[ur]^{\varphi}\ar[dr]_{\varphi^{\prime}}\\
&X^{\prime}}
\end{equation}
It is clear that the pair $(*, id)$ is the initial object in $\mathcal{C}(*)$. Therefore, the category $\mathcal{C}(*)$ is contractible by Lemma \ref{elem}(iii).

There is a natural functor $F\colon \mathcal{C}(*)\to \mathcal{P}_R\mathcal{C}$, which assigns the 1-path $*\xrightarrow{\varphi} X\xleftarrow{id}X$ to an object $*\xrightarrow{\varphi}X$ of $\mathcal{C}(*)$. This functor is the right adjoint to the functor which assigns to the path $*\xrightarrow{\varphi} X\leftarrow Y$ the object $*\xrightarrow{\varphi} X$. Therefore, the categories $\mathcal{P}_R\mathcal{C}$ and $\mathcal{C}(*)$ are homotopy equivalent by Lemma \ref{elem}(ii).
\endproof

\subsubsection{The Quillen's Theorem B}
We are going to prove that the natural projection functor $p\colon \mathcal{P}_R\EXT^k(M,N)\to\EXT^k(M,N)$
satisfies the conditions of the Quillen's Theorem B. Let us recall what it is. Actually we formulate the Corollary to Theorem B which we
will only need.

Let $g\colon E\to B$ be a map of topological spaces, $b\in B$.
The {\it homotopical fiber} $F(g,b)$ of the map $g$ over $b$ is the spaces
\begin{equation}\label{homfib}
F(g,b)=\{e\in E \text{ and a path } p\colon [0,1]\to B|\ p(0)=g(e),\ p(1)=b\}
\end{equation}
In other words, we replace a map $g\colon E\to B$ by a homotopical to it Serre fibration, and then compute the fiber.
If the map $g\colon E\to B$ is a Serre fibration itself, then up to homotopy the homotopical fiber coincides with the set-theoretical
fiber.

\begin{example}
Let $G$ be a topological group considered as a category $\underline{G}$ with a single object, and let $BG$ be its classifying space.
There is a bundle $p\colon EG\to BG$ which is a Serre fibration with fiber $G$, and with $EG$ homotopically trivial. Then $G$ is homotopically
equivalent to the homotopical fiber, which is $\Omega BG$ since $EG$ is contractible.
$$
G\simeq\Omega BG
$$
\end{example}

Let $f\colon \mathcal{C}_1\to\mathcal{C}_2$ be a functor. For any $Y\in\mathrm{Ob}(\mathcal{C}_2)$ define the
{\it comma category} $Y\setminus f$ whose objects are the pairs $\{(X,v)|X\in \mathrm{Ob}\mathcal{C}_1,\ v\in\Mor_{\mathcal{C}_2}(Y,f(X))\}$.
The morphisms are the maps $X\to X^\prime$ such that the corresponding triangle is commutative.

Morally, the comma category $Y\setminus f$ has the role of the fiber category
$f^{-1}(Y)$ consisting from objects over $Y$ and their morphisms which are mapped by $f$ to the identity of $Y$.
Let us discuss a condition under which the categories $Y\setminus f$ and $f^{-1}(Y)$ are indeed homotopically equivalent.

One says that a functor $f\colon \mathcal{C}_1\to\mathcal{C}_2$ makes the category
$\mathcal{C}_1$ a {\it prefibred category} over $\mathcal{C}_2$ if for every object $Y\in\mathrm{Ob}\mathcal{C}_2$ the functor
\begin{equation}\label{prefibred}
f^{-1}(Y)\to Y\setminus f,\ \ \ X\mapsto (X,id_Y)
\end{equation}
has a right adjoint.

It follows from Lemma \ref{elem}(ii) that the functor above is a homotopical equivalence.

In the case above denote by $(X,v)\mapsto v^*X$ the right adjoint functor to\eqref{prefibred}. Then for any arrow $u\colon Y\to Y^\prime$
in $\mathcal{C}_2$ one has the following functor $u^*\colon f^{-1}(Y^\prime)\to f^{-1}(Y)$:
\begin{equation}\label{prefibred2}
f^{-1}(Y^\prime)\rightarrow Y^\prime\setminus f \xrightarrow{u} Y\setminus f \xrightarrow{v^*} f^{-1}(Y)
\end{equation}
This functor is called {\it the base change functor}.

Now we are ready to formulate the {\bf Corollary to the Quillen's Theorem B}:
\begin{theorem}\label{quillenB}
Suppose a functor $f\colon \mathcal{C}_1\to\mathcal{C}_2$ is prefibred and that for any arrow $u\colon Y\to Y^\prime$ in $\mathcal{C}_2$
the base change functor $u^\prime\colon f^{-1}(Y^\prime)\to f^{-1}(Y)$ is a homotopy equivalence.
Then for any $Y\in\mathrm{Ob}\mathcal{C}_2$, the classifying space of the category $f^{-1}(Y)$ is homotopy equivalent to the homotopy-fibre of
$f\colon B\mathcal{C}_1\to B\mathcal{C}_2$ over $Y$.
\end{theorem}
See [Q], Section 1 for a proof.
\endproof

\subsubsection{The functor $p\colon \mathcal{P}_R\EXT^k(M,N)\to\EXT^k(M,N)$ satisfies the hypothesis of the Quillen's Theorem B}
Let $Y\in\mathrm{Ob}\EXT^k(M,N)$. Consider the category $Y\setminus p$ (see previous Subsection). There is a natural inclusion of categories
$i\colon p^{-1}(Y)\to Y\setminus p$. We firstly prove
\begin{lemma}
The functor $i$ has a right adjoint.
\end{lemma}
\proof{}
This is very easy. Let
\begin{equation}\label{bis10}
\xymatrix{
{*}\ar[r] &X^{\prime}&Y^{\prime}\ar[l]_{\varphi}\\
&&Y\ar[u]_{v}}
\end{equation}
be an element of the category $Y\setminus p$. The right adjoint functor $v^*$ associates to it the 1-path
\begin{equation}\label{bis11}
*\rightarrow X^{\prime}\xleftarrow{\varphi\circ v}Y
\end{equation}
\endproof

Now for an arrow $u\colon Y\to Y^{\prime}$ the base change functor $u^\prime\colon p^{-1}(Y^\prime)\to p^{-1}(Y)$ just maps
\begin{equation}\label{bis12}
*\rightarrow X\prime\xleftarrow{\varphi}Y^\prime
\end{equation}
to
\begin{equation}\label{bis13}
*\rightarrow X^\prime\xleftarrow{\varphi\circ u}Y
\end{equation}
We prove
\begin{lemma}
For any arrow $u\colon Y\to Y^\prime$, the base change functor $u^\prime\colon p^{-1}(Y^\prime)\to p^{-1}(Y)$ is a homotopy equivalence.
\end{lemma}
\proof{}
To prove that it is a homotopy equivalence, it is enough to prove that it has a left adjoint.
We use the following Lemma due to V.Retakh [R], Lemmata 1,2:
\begin{lemma}\label{pushouts}
\begin{itemize}
\item[1.] Every morphism in $\EXT^k(M,N)$ can be canonically factorized as the composition of two morphisms such that the first is (component-wise) monomorphic, and the second admits a section;
\item[2.] pushouts along component-wise monomorphic morphisms exist in $\EXT^k(M,N)$.
\end{itemize}
\end{lemma}
See [Sch], Lemma 4.4 for a proof.
\endproof
Now the left adjoint functor assigns to the diagram
\begin{equation}\label{bis14}
\xymatrix{
&&Y^{\prime}\\
{*}\ar[r]& X&Y\ar[l]\ar[u]}
\end{equation}
the upper line of the diagram
\begin{equation}\label{bis15}
\xymatrix{
{*}\ar[r]& {X\sqcup_YY^\prime}&Y^\prime\ar[l]\\
{*}\ar[u]^{id}\ar[r]&X\ar[u]&Y\ar[l]\ar[u]}
\end{equation}
if the morphism $Y\to Y^\prime$ is monomorphic, in this case $X\sqcup_YY^\prime$ is an extension (that is, is exact sequence) by Lemma
\ref{pushouts}(2). In the case of general morphism $Y\to Y^\prime$ we firstly use Lemma \ref{pushouts}(1) to replace \eqref{bis15} by a diagram with
monomorphic $Y\to Y^\prime$.
\endproof

We have proved that the functor $p\colon \mathcal{P}_R\EXT^k(M,N)\to\EXT^k(M,N)$ satisfies the Quillen's Theorem B.
\endproof

\subsubsection{The category $\Omega_R\EXT^k(M,N)=p^{-1}(*)$}
It remains to compute the fiber for some choice of the based object.

We choose for $k\ge 2$ the based object to be equal to
\begin{equation}\label{bis16}
0\rightarrow N\xrightarrow{id}N\rightarrow 0\rightarrow\dots\rightarrow 0\rightarrow M\xrightarrow{id}M\rightarrow 0
\end{equation}
and for $k=1$ to be equal to
\begin{equation}\label{bis17}
0\rightarrow N\rightarrow N\oplus M\rightarrow M\rightarrow 0
\end{equation}
with the natural inclusion and projection.

Denote the category $p^{-1}(*)$ by $\Omega_R\EXT^k(M,N)$.
We prove
\begin{lemma}
The category $\Omega_R\EXT^k(M,N)$ is homotopy equivalent to the category $\EXT^{k-1}(M,N)$.
\end{lemma}
\proof{}
Consider firstly the case $k\ge 2$.
Suppose $*\rightarrow X\leftarrow *$ be an object in $\Omega_R\EXT^k(M,N)$, where $*$ is given by \eqref{bis16}, and $X$ is an extension
\begin{equation}\label{bis18}
0\rightarrow N\rightarrow X_1\rightarrow\dots\rightarrow X_k\rightarrow M\rightarrow 0
\end{equation}
In particular, we have the following diagram:
\begin{equation}\label{bis19}
\xymatrix{
M\ar[r]^{id}\ar[d]^{\alpha}&M\ar[d]^{id}\ar[r]&0\\
X_k\ar[r]^{\delta_k}&M\ar[r]&0\\
M\ar[u]_{\beta}\ar[r]^{id}&M\ar[u]_{id}\ar[r]&0}
\end{equation}
In particular, both $\alpha$ and $\beta$ give sections of $\delta_k$. We have an extension and a map
\begin{equation}\label{bis20}
\xymatrix{
0\ar[r]&N\ar[r]&X_1\ar[r]&\dots\ar[r]&X_{k-1}\ar[r]&\mathrm{Ker}\delta_k\ar[r]&0\\
&&&&&M\ar[u]_{\alpha-\beta}}
\end{equation}
The pull-back of this extension by $\alpha-\beta$ gives an element in $\EXT^{k-1}(M,N)$. We have constructed a functor
$F\colon\Omega_R\EXT^k(M,N)\to\EXT^{k-1}(M,N)$. Let us prove that it is a homotopy equivalence.

By Lemma \ref{elem}(ii), it is enough to construct a left adjoint functor $G$.

Suppose
\begin{equation}\label{bis21}
0\rightarrow N\rightarrow Y_1\rightarrow\dots\rightarrow Y_{k-1}\xrightarrow{\delta_{k-1}} M\rightarrow 0
\end{equation}
be a $(k-1)$-extension.
Construct the extension
\begin{equation}\label{bis22}
0\rightarrow N\rightarrow Y_1\rightarrow\dots\rightarrow Y_{k-1}\xrightarrow{\Delta\circ\delta_{k-1}} M\oplus M\xrightarrow{(1,-1)}M\rightarrow 0
\end{equation}
where $\Delta\colon M\hookrightarrow M\oplus M$ is the diagonal imbedding. We can define by the last extension an element in $\Omega_R\EXT^k$ in a natural way. This functor $G\colon\EXT^{k-1}(M,N)\to\Omega_R\EXT^k(M,N)$ is the left adjoint to $F$.
\endproof
The first two claims of Theorem \ref{theor1} are proven. Concerning the third claim we refer the reader to [NR], Section 8.

The case $k=1$ is trivial, because there is only one loop up to an isomorphism. On the other hand, $\EXT^0(M,N)$ is a discrete category with only identity maps and the set of objects equal to $\Hom(M,N)$, by definition.
\endproof

\section{A construction of S.Schwede}
Here we describe the construction of Stefan Schwede [Sch] which
mainly motivated our work.
\subsection{}
We are going to apply the results of Section 1 in the case when the
abelian category $\mathcal{A}$ (in which we consider extensions) is
also {\it monoidal}. This means that there is a bifunctor
$\otimes\colon \mathcal{A}\times \mathcal{A}\to\mathcal{A}$ which is
associative and distributive with respect to the direct sum. We
always suppose that there is a two-sided unit object for the
monoidal structure. A priori the monoidal structure is not exact,
but we want to consider it as it would be. For this we impose the
following condition:

\begin{itemize}
\item[(*)]
There is a full additive subcategory
$\mathcal{A}_0\subset\mathcal{A}$ such that
\begin{itemize}
\item[(i)] the monoidal structure
is exact on $\mathcal{A}_0$;
\item[(ii)] the categories
$\EXT^k_{\mathcal{A}_0}$ and $\EXT^k_\mathcal{A}$ are homotopically
equivalent for any $k$ (the inclusion of the nerves is a weak homotopy equivalence);
\item[(iii)] the unit object belongs to
$\mathcal{A}_0$;
\item[(iv)] the category $\mathcal{A}_0$ is closed under the monoidal structure in $\mathcal{A}$; hence, it is itself a monoidal additive category.
\end{itemize}
\end{itemize}

\begin{example}
Let $A$ be an associative algebra, and let $\mathcal{A}$ be the
category of $A$-bimodules. The monoidal structure is the tensor
product of bimodules over $A$, $(M_1,M_2)\mapsto M_1\otimes_A M_2$.
This product is not exact. Let $\mathcal{A}_0$ be the full additive
subcategory of bimodules which are flat as right $A$-modules. The condition (*) for $\mathcal{A}_0$ is proven in
[Sch], Lemma 2.1. The unit object is the tautological bimodule $A$, and it
belongs to $\mathcal{A}_0$. It is clear that such $\mathcal{A}_0$ is closed under the monoidal structure.
\end{example}

\begin{example}
Suppose $A$ is an associative bialgebra (see definition below in
Section 4.1), and let $\mathcal{A}$ be the category of left $A$-modules (as
over algebra). For any two left modules $M_1$ and $M_2$, their
tensor product $M_1\otimes_k M_2$ over the ground field is naturally
a module over the algebra $A\otimes_k A$. Now the coproduct
$\Delta\colon A\to A\otimes_k A$ (which is a map of algebras) makes
$M_1\otimes M_2$ a left $A$-module. This monoidal structure is
exact, one can set $\mathcal{A}_0=\mathcal{A}$. The unit object is
the trivial module $k$; for existence of it one needs to have the
counit.
\end{example}

The two ingredients of our game are the Yoneda product and the
Schwede's tensor product of extensions.

\subsubsection{The Yoneda product}
For extensions $\mathcal{E}=\{0\rightarrow M\rightarrow
E_1\rightarrow\dots E_k\rightarrow N\rightarrow 0\}$ and
$\mathcal{F}=\{0\rightarrow N\rightarrow
F_1\rightarrow\dots\rightarrow F_\ell\rightarrow P\rightarrow 0\}$
their Yoneda product $\mathcal{E}\sharp\mathcal{F}$ is the extension
\begin{equation}\label{eq2.3}
0\rightarrow M\rightarrow E_1\rightarrow\dots E_k\rightarrow
F_1\rightarrow\dots\rightarrow F_\ell\rightarrow P\rightarrow 0
\end{equation}
where the ``central'' arrow $E_k\to F_1$ is the composition
$E_k\rightarrow N\rightarrow F_1$. Clearly it is an extension.

\begin{lemma}
Under the map $\varphi$ of Lemma 1.1, the Yoneda product is
corresponded to the natural product $\Ext^k(N,M)\otimes
\Ext^\ell(P,N)\to \Ext^{k+\ell}(P,M)$ (which also will be called the
Yoneda product).
\end{lemma}
\endproof
\subsubsection{The Schwede's tensor product}
Let now $\mathcal{E}$ and $\mathcal{F}$ be extensions as above, with
$M=N=P$ equal to $A$, the unit object of the monoidal category
$\mathcal{A}$. Suppose also that the terms of these extensions
belong to the subcategory $\mathcal{A}_0$, on which $\otimes$ is
exact.

The idea is to define an extension
$\mathcal{E}{\otimes}_\tau\mathcal{F}$ of length $k+\ell$, such that
there is a diagram of maps of extensions:
\begin{equation}\label{eq2.4}
\xymatrix{ & \mathcal{E}{\otimes}_\tau\mathcal{F}\ar[dl]\ar[dr]\\
\mathcal{E}\sharp\mathcal{F}&&(-1)^{k\ell}\mathcal{F}\sharp\mathcal{E}}
\end{equation}
Then it will follow from Lemma 1.1 that the extensions
$\mathcal{E}\sharp\mathcal{F}$ and
$(-1)^{k\ell}\mathcal{F}\sharp\mathcal{E}$ define the same element
in $\Ext_\mathcal{A}^\mb(A,A)$ (in other words, the Yoneda product
in $\Ext^\mb(A,A)$ is graded commutative).

The most naive candidate for $\mathcal{E}{\otimes}_\tau\mathcal{F}$,
the usual tensor product of the complexes $\mathcal{E}$ and
$\mathcal{F}$, has the length for one more than
$\mathcal{E}\sharp\mathcal{F}$ and $\mathcal{F}\sharp\mathcal{E}$.
We modify this naive definition, as follows.

Denote by $\tau(\mathcal{E})$ the following ``truncated'' complex:
$\tau(\mathcal{E})=\{0\rightarrow A\rightarrow E_1\rightarrow\dots
E_k\rightarrow  0\}$ (what is truncated is the last term), and
analogously $\tau(\mathcal{F})=\{0\rightarrow A\rightarrow
F_1\rightarrow\dots\rightarrow F_\ell\rightarrow 0\}$. We have:
$\tau(\mathcal{E})$ is quasiisomorphic to $A[-k]$, and
$\tau(\mathcal{F})$ is quasiisomorphic to $A[-\ell]$. As the lements
of these complexes belong to $\mathcal{A}_0$, the usual tensor
product $\tau(\mathcal{E})\otimes\tau(\mathcal{F})$ is
quasi-isomorphic to $A[-k-\ell]$. This gives us the following
extension of length $k+\ell$: $0\rightarrow
\tau(\mathcal{E})\otimes\tau(\mathcal{F})\rightarrow A\rightarrow
0$. We denote this extension by
$\mathcal{E}\otimes_\tau\mathcal{F}$.
\begin{lemma}
There exists a diagram of extensions \eqref{eq2.4}.
\end{lemma}
\proof{} Represent $\tau(\mathcal{E})\otimes\tau(\mathcal{F})$ as a
rectangle, see Figure 1. The terms written down in the two marked on
the Figure borders are the quotient-complexes isomorphic to
$\mathcal{E}\sharp\mathcal{F}$ and
$(-1)^{k\ell}\mathcal{F}\sharp\mathcal{E}$, correspondingly.
\sevafigc{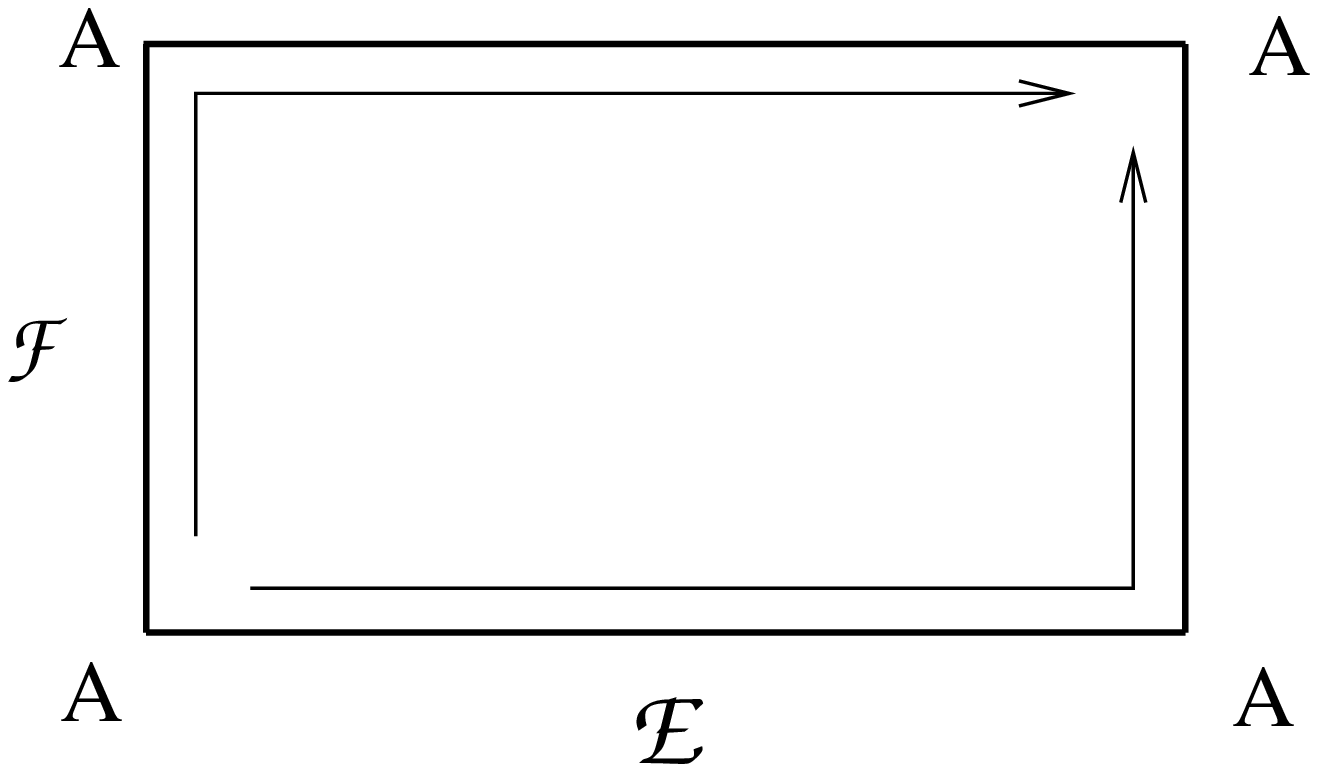}{60mm}{0}{}
\endproof

\begin{corr}
Let $\mathcal{A}$ be a monoidal category satisfying the condition
(*), and let $A$ be the unit object in it. Then the Yoneda product
on $\Ext^\mb_\mathcal{A}(A,A)$ is graded commutative.
\end{corr}
It follows from the Lemma above, Lemma 2.1, and Lemma 1.1.
\endproof
\subsection{The Lie bracket}
One can go one step further in these constructions, and consider the
following (non-commutative in any sense) diagram:

\begin{equation}\label{eq2.5}
\xymatrix{ & \mathcal{E}{\otimes}_\tau\mathcal{F}\ar[dl]\ar[dr]\\
\mathcal{E}\sharp\mathcal{F}&&(-1)^{k\ell}\mathcal{F}\sharp\mathcal{E}\\
&(-1)^{k\ell}\mathcal{F}\otimes_\tau\mathcal{E}\ar[ul]\ar[ur]}
\end{equation}

This diagram defines a loop in $\mathcal{N}\EXT^{k+\ell}(A,A)$, that
is, an element in $\pi_1(\mathcal{N}\EXT^{k+\ell}(A,A)$. By Retakh's
Theorem 1.3 this $\pi_1$ is {\it canonically} isomorphic to
$\Ext^{k+\ell-1}_\mathcal{A}(A,A)$. What we get is the bracket,
induced by the Gerstenhaber bracket on Hochschild cohomology:

\begin{theorem}(Schwede)
The Lie bracket on $HH^\mb(A,A)=\Ext^\mb_{A-\Bimod}(A,A)$ defined
above is equal to the bracket induced from the Gerstenhaber bracket.
In particular, it obeys the graded Jacobi identity, and together
with the Yoneda product forms a 2-algebra.
\end{theorem}

\begin{remark}
In the above construction, one place should be treated more
carefully. Namely, the Retakh's theorem was proven for some fixed
choice of the basepoint. The categories $\EXT^k(A,A)$ are not
connected ($\pi_0(\EXT^k(A,A))=\Ext^k(A,A)$), and the computation of
the homotopy groups may depend on the choice of connection
component. For this we use Lemma \ref{allequivalent} which shows that all
connection components are homotopically equivalent.
\end{remark}

\section{$n$-fold monoidal categories}
\subsection{Introduction}
An $n$-fold monoidal category $\mathcal{C}$ is a category with $n$
monoidal structures $\otimes_1,\dots,\otimes_n\colon
\mathcal{C}\times\mathcal{C}\to\mathcal{C}$ which obey some
compatibility relations. They were introduced in [BFSV] with the
following motivation.

The classifying space of a 1-monoidal category admits an action of
the Stasheff operad (see [M2]). A {\it connected} space with an
action of this operad has a homotopical type of loop spaces; the
same is true for non-connected spaces whose linear connection
components form a group (see loc.cit.). In [BFSV], the authors found
a structure on a category (called $n$-monoidal), such that there is
an operad of categories acting on it; the classifying space of this
operad of categories is $n$-dimensional little disc operad. This is
very close to say that the classifying space of an $n$-monoidal
category is an $n$-fold loop space (see [M2]). The definition is
iterative, like the definition of $n$-fold loop space.

\subsection{Definition}
\begin{defn}
A (strict) monoidal category is a category $\mathcal{C}$ together
with a functor $\otimes\colon
\mathcal{C}\times\mathcal{C}\to\mathcal{C}$ and an object
$\mho\in\Ob(\mathcal{C})$ such that
\begin{itemize}
\item[1.] $\otimes$ is strictly associative;
\item[2.] $\mho$ is a strict two-sided unit for $\otimes$.
\end{itemize}
A monoidal functor $(F,\eta)\colon\mathcal{C}\to\mathcal{D}$ between
monoidal categories is a functor $F$ such that
$F(\mho_\mathcal{C})=\mho_\mathcal{D}$ with a natural transformation
\begin{equation}\label{eta}
\eta_{A,B}\colon F(A)\otimes F(B)\to F(A\otimes B)
\end{equation}
which satisfies the following conditions:
\begin{itemize}
\item[1.] Internal associativity: the following diagram commutes \begin{equation}\label{intassoc}
\xymatrix{ F(A)\otimes F(B)\otimes F(C)\ar[rr]^{\eta_{A,B}\otimes
id_{F(C)}}\ar[d]^{id_{F(A)}\otimes\eta_{B,C}}&&F(A\otimes B)\otimes
F(C)\ar[d]^{\eta_{A\otimes B,C}}\\
F(A)\otimes F(B\otimes C)\ar[rr]^{\eta_{A,B\otimes C}}&&F(A\otimes
B\otimes C) }
\end{equation}
\item[2.] Internal unit conditions:
$\eta_{A,\mho}=\eta_{\mho,A}=id_{F(A)}$.
\end{itemize}
\end{defn}
The crucial in this definition is that the map $\eta$ is not
required to be an isomorphism.

Denote by $\mathbf{MonCat}$ the category of (small) monoidal
categories and monoidal functors.
\subsubsection{}
\begin{defn}
A 2-fold monoidal category is a monoid in $\mathbf{MonCat}$. This
means, that we are given a monoidal category
$(\mathcal{C},\otimes_1,\mho)$, and a monoidal functor
$(\otimes_2,\eta)\colon\mathcal{C}\times\mathcal{C}\to\mathcal{C}$
which satisfies the following axioms:
\begin{itemize}
\item[1.] External associativity: the following diagram commutes in
$\mathbf{MonCat}$
\begin{equation}\label{external}
\xymatrix{
\mathcal{C}\times\mathcal{C}\times\mathcal{C}\ar[rr]^{(\otimes_2,\eta)\times
id_{\mathcal{C}}}\ar[d]^{id_{\mathcal{C}}\times(\otimes_2,\eta)}&&\mathcal{C}\times\mathcal{C}\ar[d]^{(\otimes_2,\eta)}\\
\mathcal{C}\times\mathcal{C}\ar[rr]^{(\otimes_2,\eta)}&&\mathcal{C}}
\end{equation}
\item[2.] External unit conditions:
the following diagram commutes in $\mathbf{MonCat}$
\begin{equation}\label{externalunit}
\xymatrix{
\mathcal{C}\times\mho\ar[r]^{\subseteq}\ar[d]^{\cong}&\mathcal{C}\times\mathcal{C}\ar[d]^{(\otimes_2,\eta)}&\mho\times\mathcal{C}
\ar[l]_{\supseteq}\ar[d]^{\cong}\\
\mathcal{C}\ar[r]^{=}&\mathcal{C}&\mathcal{C}\ar[l]_{=}}
\end{equation}
\end{itemize}
\end{defn}
Let us note that the role of the monoidal structures $\otimes_1$ and
$\otimes_2$ in this definition is not symmetric.

Explicitly the definition above means that we have an operation
$\otimes_2$ with the two-sided unit $\mho$ (the same that for
$\otimes_1$) and a natural transformation
\begin{equation}\label{natural}
\eta_{A,B,C,D}\colon (A\otimes_2 B)\otimes_1(C\otimes_2 D)\to
(A\otimes_1 C)\otimes_2 (B\otimes_1 D)
\end{equation}
The internal unit conditions are:
$\eta_{A,B,\mho,\mho}=\eta_{\mho,\mho,A,B}=id_{A\otimes_2B}$, and
the external unit conditions are:
$\eta_{A,\mho,B,\mho}=\eta_{\mho,A,\mho,B}=id_{A\otimes_1B}$. As
well, one has the morphisms
\begin{equation}\label{12left}
\eta_{A,\mho,\mho,B}\colon A\otimes_1B\to A\otimes_2B
\end{equation}
and
\begin{equation}\label{12right}
\eta_{\mho,A,B,\mho}\colon A\otimes_1B\to B\otimes_2A
\end{equation}
which will be very essential in Section 3.

The internal associativity gives the commutative diagram:
\begin{equation}\label{intexpl}
\xymatrix{
(U\otimes_2V)\otimes_1(W\otimes_2X)\otimes_1(Y\otimes_2Z)\ar[rrr]^{\eta_{U,V,W,X}\otimes_1id_{Y\otimes_2Z}}
\ar[dd]^{id_{U\otimes_2V}\otimes_1\eta_{W,X,Y,Z}}&&&
\bigl((U\otimes_1W)\otimes_2(V\otimes_1X)\bigr)\otimes_1(Y\otimes_2Z)\ar[dd]^{\eta_{U\otimes_1W,V\otimes_1X,Y,Z}}\\
\\
(U\otimes_2V)\otimes_1\bigl((W\otimes_1Y)\otimes_2(X\otimes_1Z)\bigr)\ar[rrr]^{\eta_{U,V,W\otimes_1Y,X\otimes_1Z}}&&&
(U\otimes_1W\otimes_1Y)\otimes_2(V\otimes_1X\otimes_1Z)}
\end{equation}
The external associativity condition gives the commutative diagram:
\begin{equation}\label{extexpl}
\xymatrix{(U\otimes_2V\otimes_2W)\otimes_1(X\otimes_2Y\otimes_2Z)\ar[rrr]^{\eta_{U\otimes_2V,W,X\otimes_2Y,Z}}
\ar[dd]^{\eta_{U,V\otimes_2W,X,Y\otimes_2Z}}&&&\bigl((U\otimes_2V)\otimes_1(X\otimes_2Y)\bigr)\otimes_2(W\otimes_1Z)
\ar[dd]^{\eta_{U,V,X,Y}\otimes_2id_{W\otimes_1Z}}\\
\\
(U\otimes_1X)\otimes_2\bigl((V\otimes_2W)\otimes_1(Y\otimes_2Z)\bigr)\ar[rrr]^{id_{U\otimes_1X}\otimes_2\eta_{V,W,Y,Z}}&&&
(U\otimes_1X)\otimes_2(V\otimes_1Y)\otimes_2(W\otimes_1Z) }
\end{equation}

Finally, [BFSV] gives
\begin{defn}
Denote by $\mathbf{MonCat}_n$ the category of (small) $n$-fold
monoidal categories. Then an $(n+1)$-fold monoidal category is a
monoid in $\mathbf{MonCat}_n$.
\end{defn}

This gives the following compatibility axiom: for $1\le i<j<k\le n$ the following diagram is commutative:

\begin{equation}\label{compexpl}
{\scriptsize\xymatrix{
&\bigl(((A_1\otimes_kA_2)\otimes_j(B_1\otimes_kB_2)\bigr)\otimes_i\bigl((C_1\otimes_kC_2)\otimes_j(D_1\otimes_kD_2)\bigr)
\ar[ldd]^{\eta^{jk}\otimes_i\eta^{jk}}\ar[dd]_{\eta^{ij}}\\
\\
\bigl((A_1\otimes_jB_1)\otimes_k(A_2\otimes_jB_2)\bigr)\otimes_i\bigl((C_1\otimes_jD_1)\otimes_k(C_2\otimes_j D_2)\bigr)\ar[dd]^{\eta^{ik}}&
\bigl((A_1\otimes_k A_2)\otimes_i(C_1\otimes_kC_2)\bigr)\otimes_j\bigl((B_1\otimes_kB_2)\otimes_i(D_1\otimes_kD_2)\bigr)\ar[dd]_{\eta^{ik}\otimes_j\eta^{ik}}\\
\\
\bigl((A_1\otimes_jB_1)\otimes_i(C_1\otimes_jD_1)\bigr)\otimes_k\bigl((A_2\otimes_jB_2)\otimes_i(C_2\otimes_jD_2)\bigr)\ar[ddr]^{\eta^{ij}\otimes_k\eta^{ij}}&
\bigl((A_1\otimes_iC_1)\otimes_k(A_2\otimes_iC_2)\bigr)\otimes_j\bigl((B_1\otimes_iD_1)\otimes_k(B_2\otimes_iD_2)\bigr)\ar[dd]_{\eta^{jk}}\\
\\
&\bigl((A_1\otimes_iC_1)\otimes_j(B_1\otimes_iD_1)\bigr)\otimes_k\bigl((A_2\otimes_iC_2)\otimes_j(B_2\otimes_iD_2)\bigr)
}}
\end{equation}

\subsection{Examples}
See examples of monoidal (=1-monoidal) categories given in Examples
2.1 and 2.2.

\begin{example}
According to a result of Joyal and Street [JS], there are just few examples for $n\ge 2$ when the map $\eta^{ij}_{A,B,C,D}$ are {\it isomorphisms}
for any $A,B,C,D$ and any $1\le i<j\le n$ and when {\it there is a common unit object for all $n$ monoidal structures}. For $n=2$ any such category is equivalent as a 2-fold monoidal category to a category
with $A\otimes_1 B=A\otimes_2 B$ with a {\it braiding} $c_{A,B}\colon A\otimes B\to B\otimes A$ defining a structure of a {\it braided category} (see, e.g., [ES]) on
$\mathcal{C}$. Then we can construct a map $\eta_{A,B,C,D}\colon
(A\otimes B)\otimes(C\otimes D)\to (A\otimes C)\otimes (B\otimes D)$
just as $\eta_{A,B,C,D}=id\otimes c_{23}\otimes id$. This
construction gives a 2-fold monoidal category. For $n>2$ and $\eta_{A,B,C,D}$ isomorphisms one necessarily
has $A\otimes_i B=A\otimes_j B$ for any $i,j$ and {\it all $\otimes_i$ are symmetric}.
We recommend the Fiedorowicz's Obervolfach talk [F] (page 4 and thereafter) for a concise but clear overview of this result.
\end{example}

\begin{example}\label{extetra}
Let $A$ be an associative bialgebra. We define a {\it tetramodule}
over it as a $k$-vector space $M$ such that there is a bialgebra
structure on $A\oplus\epsilon M$, where $\epsilon^2=0$ and the
restriction of the bialgebra structure to $A$ is the initial one
(see Section 4 for details). If we perform this definition replacing
``bialgebra'' by ``associative algebra'', we recover the concept of
bimodule; thus, this Example is a generalization of Example 2.1. We
construct in Section 4 a 2-fold monoidal structure on the abelian
category $\Tetra(A)$ of tetramodules over $A$.
\end{example}

\begin{example}\label{monalg}
Examples 2.1 and \ref{extetra} can be generalized as follows. Recall
from Example 2.2 that the left modules over an associative bialgebra
form a monoidal category, with the monoidal structure equal to the
tensor product of the underlying vector spaces. Define an {\it
$n$-fold monoidal bialgebra} as an associative algebra with $n$
coassociative coproducts $\Delta_1,\dots,\Delta_n\colon A\to
A\otimes_k A$ such that the corresponding $n$ monoidal structures on
the category of left $A$-modules form an $n$-fold monoidal category.
Thus, 0-monoidal bialgebra is just an associative algebra, and
1-monoidal bialgebra is a bialgebra. One can define the category of
tetramodules over an $n$-monoidal bialgebra analogously to the
previous Example. We claim that this category is an $(n+1)$-fold
monoidal category; a proof will appear somewhere. The author
believes that this concept of $n$-monoidal bialgebra is a conceptually right
$n$-categorical generalization of the concept of bialgebra.
\end{example}

\subsection{The operad of categories governing the $n$-fold monoidal
categories}
Fix $n\ge 1$. For any $d\ge 0$ denote by $\mathcal{M}_n(d)$ the full subcategory of the free $n$-fold monoidal category generated by objects $x_1,\dots, x_d$ consisting of objects which are monomials in $x_i$, where each $x_i$ occurs exactly ones.
For example, such monomials for $d=3$ and $n=2$ could be $(x_3\otimes_1x_1)\otimes_2 x_2$, or $(x_2\otimes_2 x_3)\otimes_1 x_1$. For fixed $n$ and $d$ the category $\mathcal{M}_n(d)$ has a finite number of objects. The morphisms in $\mathcal{M}_n(d)$ are exactly those which can be obtained as compositions of the associativities for a fixed $\otimes_i$, and $\eta_{ijkl}$, with exactly the same commutative diagrams as in $n$-fold monoidal category.

When $n$ is fixed and $d$ is varied, the categories $\mathcal{M}_n(d)$ form an operad of categories.
The following lemma follows from the definitions.

\begin{lemma}
A category is $n$-fold monoidal if and only if there is an action of the operad $\{\mathcal{M}_n(d)\}_{d\ge 0}$ of categories on it.
\end{lemma}
\endproof

The following very deep theorem is in a sense the main result in [BFSV]:
\begin{theorem}
The classifying space of the operad of categories $\{\mathcal{M}_n(d)\}$ is an operad of topological space which is homotopically equivalent (as operad) to the $n$-dimensional little discs operad.
\end{theorem}

\subsection{The category of extensions of abelian $n$-fold monoidal
category} Let $\mathcal{C}$ be an abelian $n$-fold monoidal category
with common unit object $A$. Consider the category of extensions
$\EXT^k_\mathcal{C}(A,A)$. The following Lemma is in a sense one of
the main our observations.
\begin{klemma}\label{extmon}
Under the above conditions, and the condition (*) in Section 2.1,
the disjoint union of categories $\coprod_{k\ge
1}\EXT^k_\mathcal{C}(A,A)$ is an $(n+1)$-monoidal category.
\end{klemma}
\proof{} Let $\otimes_1,\dots,\otimes_n$ be the (ordered) set of $n$
monoidal structures in $\mathcal{C}$. Each of them, $\otimes_s$,
defines a monoidal structure $\otimes_{s,\tau}$ (the Schwede's
tensor product), as is explained in Section 2.1.2. Set
$\bigotimes_i=\otimes_{i,\tau}$, $i=1\dots n$, and let
$\bigotimes_{n+1}$ be the Yoneda product of extensions (see
\eqref{eq2.3}). The common identity object is the distinguished
object in $\EXT^1_\mathcal{C}(A,A)$ (that is, the extension
$0\rightarrow A\xrightarrow{id}A\rightarrow 0$). One needs to
construct the maps $\eta_{M,N,P,Q}$ for the pairs
$(\bigotimes_i,\bigotimes_{n+1})$, $1\le i\le n$, of monoidal
structures. That is, we need to construct a map
\begin{equation}\label{rectangle}
\eta_{M,N,P,Q}\colon (M\sharp N)\otimes_{i,\tau}(P\sharp Q)\to
(M\otimes_{i,\tau} P)\sharp (N\otimes_{i,\tau}Q)
\end{equation}
The idea is shown in Figure 2. \sevafigc{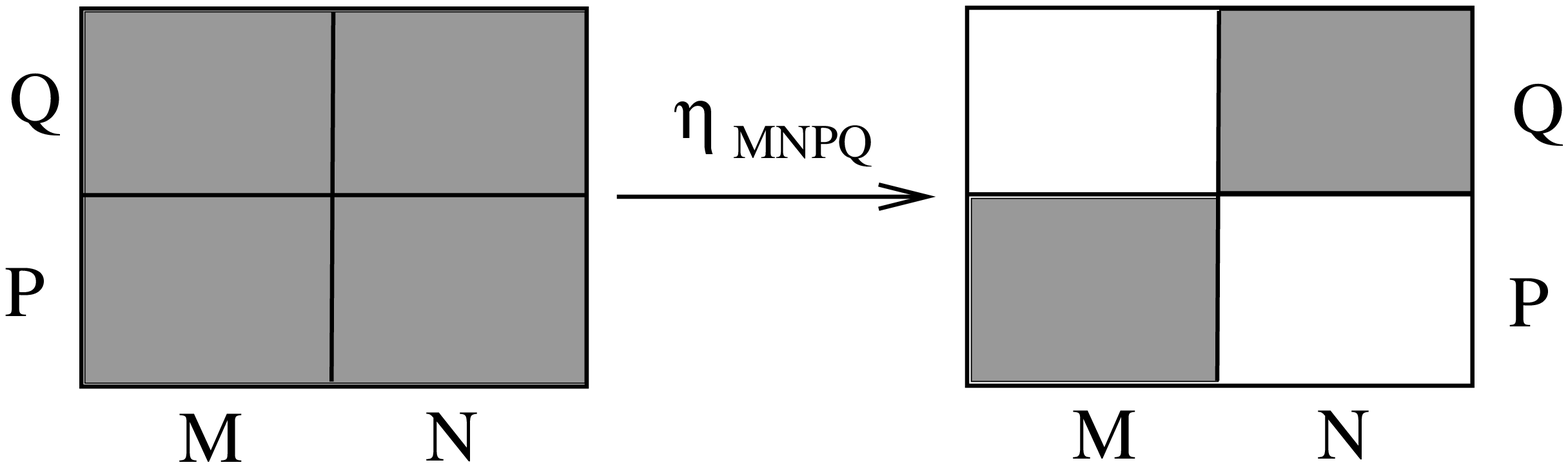}{120mm}{0}{} The
idea is to project the whole rectangle (the l.h.s. of
\eqref{rectangle}) to the marked rectangles (the r.h.s. of
\eqref{rectangle}). Precisely, this map of complexes can be
constructed as follows. Let $M$ be an $m$-extension, $N$ be an
$n$-extension, etc. Let $X(a)\otimes Y(b)$ be an element of the
rectangle, where $X$ is either $M$ or $N$, $Y$ is either
$P$ or $Q$. Note that $M(1)=P(1)=A$, while $N(1)=N_1$ and $Q(1)=Q_1$.
On the other hand, $M(m+1)=M_m$, $P(p+1)=P_p$, while $N(n+1)=Q(q+1)=A$.
Denote also by $i_M\colon A\to
M_1$ the first arrow, and by $p_M\colon M_m\to A$ the last arrow,
analogously for $N,P$ and $Q$.

One should be especially careful with the elements in $(M\sharp
N)\otimes_{i,\tau}(P\sharp Q)$ shown by the bold points in Figure 3.
The map $\eta_{M,N,P,Q}$ is zero on the two white rectangles except
these bold points.

\sevafigc{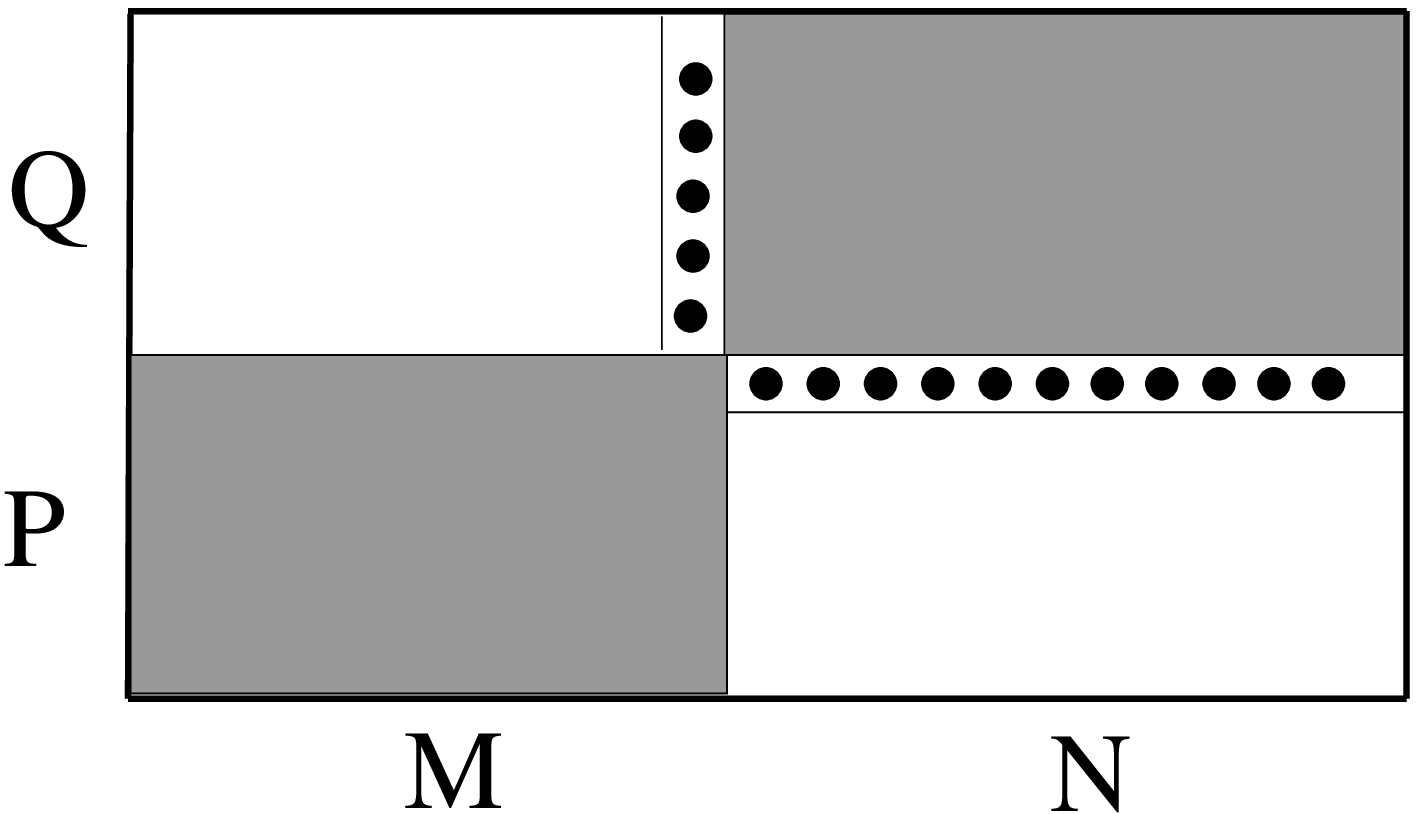}{60mm}{0}{}

\begin{remark}
In what follows we use two different indexations. The lower index $M_i$, $1\le i\le m+2$ denotes the $i$-th term of the extension $M$,
and $M(i)$, $1\le i\le m+1$ denotes a coordinate of a point of the rectangle.
\end{remark}

The map $\eta_{M,N,P,Q}$ is now defined, up to the total sign
$(-1)^{np}$, as follows:
\begin{equation}\label{rectangle2}
\begin{aligned}
\ &\bullet M(a)\otimes P(b)\mapsto M(a)\otimes P(b)\\
 &\text{ for } 1\le a\le m+1,\
1\le b\le p+1\\
\ &\bullet N(c)\otimes Q(d)\mapsto N(c)\otimes Q(d)\\
&\text{ for } 1\le c\le n,\ 1\le d\le q\\
&\bullet N(n+1)\otimes Q(d)\mapsto 0 \text{ unless }d=q, q+1\\
&\bullet N(c)\otimes Q(q+1)\mapsto 0 \text{ unless }c=n, n+1\\
&\bullet N(n+1)\otimes Q(q)\xrightarrow{id\otimes p_Q} A\otimes A\rightarrow A\\
&\bullet N(n)\otimes Q(q+1)\xrightarrow{p_N\otimes id} A\otimes A\rightarrow A\\
&\bullet N(n+1)\otimes Q(q+1)\xrightarrow{id\otimes id}A\otimes A\rightarrow A\\
&\bullet M(a)\otimes Q(d)\mapsto 0\text{ unless } a=m+1\\
&\bullet N(c)\otimes P(b)\mapsto 0\text{ unless } b=p+1\\
&\bullet M(m+1)\otimes Q(d)\xrightarrow{p_M\otimes id}A\otimes
Q(d)=N_0\otimes Q(d) \text{ for any }d\\
&\bullet N(c)\otimes P(p+1)\xrightarrow{id\otimes p_P} N(c)\otimes A=N(c)\otimes Q_0
\end{aligned}
\end{equation}
\begin{remark}
The last two lines of \eqref{rectangle2} are corresponded to the bold points on Figure 3.
\end{remark}

\begin{lemma}
The map $\eta_{M,N,P,Q}$ defined as in \eqref{rectangle2} is a map of extensions.
\end{lemma}
\proof{}
The statement of Lemma means that the map $\eta_{M,N,P,Q}$ is a map of complexes which is $id_A$ on the both ends. This is clear.
\endproof
Now we check that taking $\sharp$ as
$\otimes_{n+1}$ as above, $\eta_{M,N,P,Q}$ indeed satisfies the
axioms \eqref{intexpl}, \eqref{extexpl} and \eqref{compexpl}.

\subsubsection{Check of \eqref{extexpl}}
We only need to check the commutativity of the diagram when the second monoidal structure $\otimes_2$ is the last one, the Yoneda product $\sharp$.
(If the both monoidal structures have numbers from 1 to $n$, the commutativity follows automatically from the assumption that $\mathcal{C}$ is an $n$-monoidal category). Consider the case when $\otimes_1$ has any number from 1 to $n$, and $\otimes_2=\sharp$.

The upper-right path of the diagram \eqref{extexpl} is schematically drawn in Figure 4 below. The bold points are corresponded to the last two lines in the definition \eqref{rectangle2}. Each of the two arrow ``projects'' the gray-color rectangle into two smaller gray rectangles. The bold points are the places outside the gray-color rectangles on which the ``projection'' is not zero.
\sevafigc{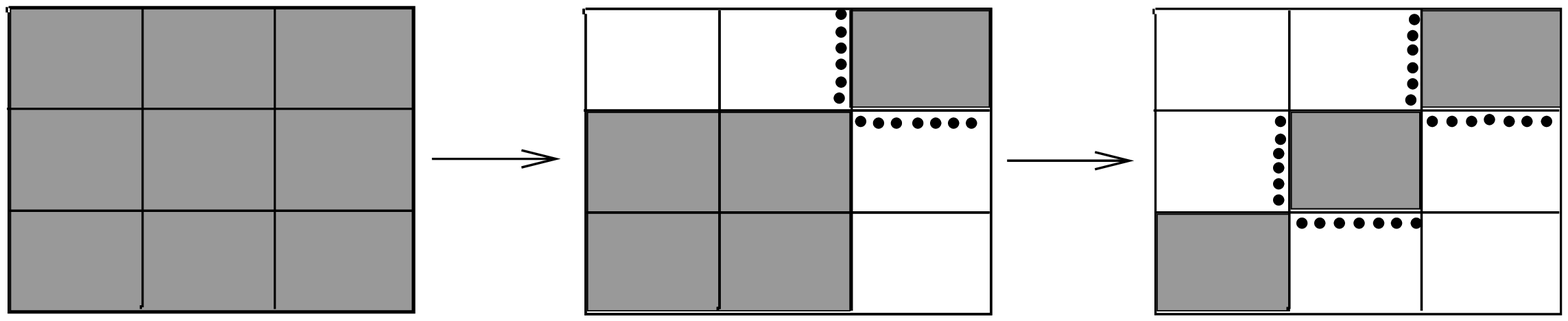}{100mm}{0}{}

The lower-left path of the diagram \eqref{rectangle2} is shown in Figure 5.
\sevafigc{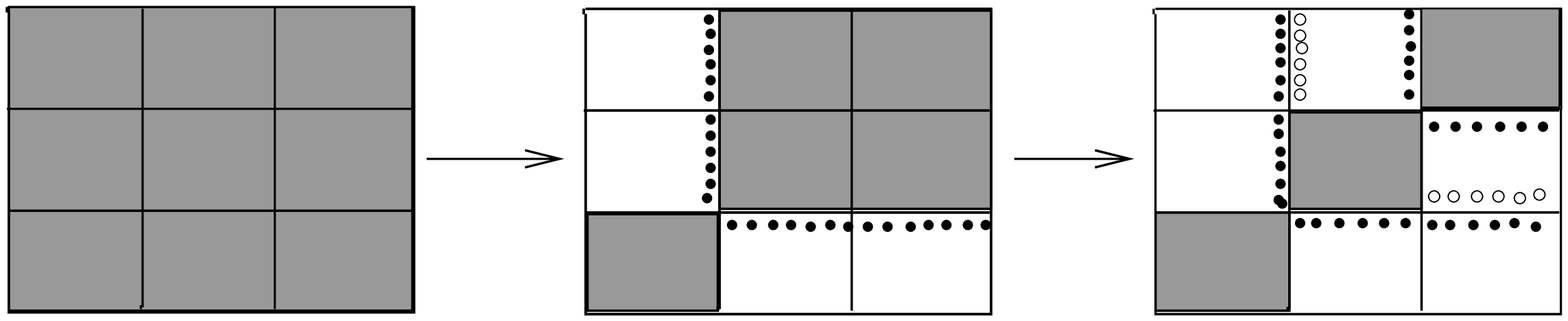}{100mm}{0}{}
After the first map, we have some bold points which we did not have in Figure 4. However, they are mapped to the corresponding white points (see Figure 5), and therefore this difference is mapped to 0 under the second arrow.

We proved, that the diagram \eqref{extexpl} is commutative in $\coprod_{k\ge
1}\EXT^k_\mathcal{C}(A,A)$.

\subsubsection{Check of \eqref{intexpl}}
The check of \eqref{intexpl} is analogous to the check of \eqref{extexpl} above, where we should consider a 2$\times$2 3-dimensional cube instead of 3$\times$3 2-dimensional cube in \eqref{extexpl}.

\subsubsection{Check of \eqref{compexpl}}
The check of \eqref{compexpl} is a bit more tricky. At first, we suppose, as above, that $k=\sharp$. Here $i$ and $j$, $i<j$, may be arbitrary.
We use the notation
\begin{equation}\label{keyvanishingelement}
(A_\alpha(\ell_1)\otimes_j B_\beta(\ell_2))\otimes_i (C_\gamma(\ell_3)\otimes_j D_\delta(\ell_4))
\end{equation}
where $\alpha,\beta,\gamma,\delta\in \{1,2\}$ for a general element of the top of the diagram \eqref{compexpl},
$((A_1\sharp A_2)\otimes_j(B_1\sharp B_2))\otimes_i ((C_1\sharp C_2)\otimes_j (D_1\sharp D_2))$.

Introduce some terminology. We say that an element \eqref{keyvanishingelement} is {\it singular}, if at least one of the maps in the (left-hand path of the) diagram belongs to the last two lines of \eqref{rectangle2} (the ``bold points''), otherwise we call this element {\it regular}.
We prove the following
\begin{lemma}\label{spisok}
Consider an element \eqref{keyvanishingelement} on which the right-hand path of \eqref{compexpl} is non-zero.
If this element is regular, it is necessarily of the form
\begin{equation}
(A_\alpha(\ell_1)\otimes_j B_\alpha(\ell_2))\otimes_i(C_\alpha(\ell_3)\otimes_jD_\alpha(\ell_4))
\end{equation}
for $\alpha\in\{1,2\}$. If this element is singular, it is necessarily of the form \eqref{keyvanishingelement} such that some number $s\in\{1,2,3\}$ among the factors $A,B,C,D$ have the corresponding $\alpha$ equal to 1, and the corresponding value of $\ell_t$ at these factors is maximal possible (we denote it by $\max$); the remaining $4-s$ factors have necessarily the lower index $\beta=2$, and their $\ell_q$ may be arbitrary. More precisely, these are elements of the following 3 groups:
\begin{equation}\label{keyvanishingfinal1}
\begin{aligned}
\ &(A_1(\max)\otimes_j B_2(\ell_1))\otimes_i(C_2(\ell_2)\otimes_j D_2(\ell_3))\\
&(A_2(\ell_1)\otimes_j B_1(\max))\otimes_i(C_2(\ell_2)\otimes_j D_2(\ell_3))\\
&(A_2(\ell_1)\otimes_j B_2(\ell_2))\otimes_i(C_1(\max)\otimes_j D_2(\ell_3))\\
&(A_2(\ell_1)\otimes_j B_2(\ell_2))\otimes_i(C_2(\ell_3)\otimes_j D_1(\max))
\end{aligned}
\end{equation}

\begin{equation}\label{keyvanishingfinal2}
\begin{aligned}
\ &(A_1(\max)\otimes_j B_1(\max))\otimes_i(C_2(\ell_1)\otimes_j D_2(\ell_2))\\
&(A_1(\max)\otimes_j B_2(\ell_1))\otimes_i(C_1(\max)\otimes_j D_2(\ell_2))\\
&(A_1(\max)\otimes_j B_2(\ell_1))\otimes_i(C_2(\ell_2)\otimes_j D_1(\max))\\
&(A_2(\ell_1)\otimes_j B_1(\max))\otimes_i(C_1(\max)\otimes_j D_2(\ell_2))\\
&(A_2(\ell_1)\otimes_j B_1(\max))\otimes_i(C_2(\ell_2)\otimes_j D_1(\max))\\
&(A_2(\ell_1)\otimes_j B_2(\ell_2))\otimes_i(C_1(\max)\otimes_j D_1(\max))
\end{aligned}
\end{equation}

\begin{equation}\label{keyvanishingfinal3}
\begin{aligned}
\ &(A_2(\ell)\otimes_j B_1(\max))\otimes_i(C_1(\max)\otimes_j D_1(\max))\\
&(A_1(\max)\otimes_j B_2(\ell))\otimes_i (C_1(\max)\otimes_j D_1(\max))\\
&(A_1(\max)\otimes_j B_1(\max))\otimes_i (C_2(\ell)\otimes_j D_1(\max))\\
&(A_1(\max)\otimes_j B_1(\max))\otimes_i (C_1(\max)\otimes_j D_2(\ell))
\end{aligned}
\end{equation}
The same statement is true for the left-hand path of \eqref{compexpl}.
\end{lemma}

\proof{}
Consider the right-hand path of diagram \eqref{compexpl}. We want to find all possible elements in the form $(A_\alpha(\ell_1)\otimes_j B_\beta(\ell_2))\otimes_i (C_\gamma(\ell_3)\otimes_j D_\delta(\ell_4))$ which are mapped to non-zero elements by the right-hand path.

The right-hand path is the composition of 3 maps. The first map does not use \eqref{rectangle2}, and, therefore, does not give any restrictions.

Consider the second map. Only the following expressions may map to nonzero elements.

{\bf Regular elements non-vanishing by the second map:}

\begin{equation}\label{keyvanishing1}
(A_\alpha(\ell_1)\otimes_i C_\alpha(\ell_2))\otimes_j(B_\beta(\ell_3)\otimes_iD_\beta(\ell_4))
\end{equation}
Here $\alpha,\beta\in\{1,2\}$, and $\ell_s$ are arbitrary.

{\bf Singular elements non-vanishing by the second map:}

\begin{equation}\label{keyvanishing2}
\begin{aligned}
\ &(A_2(\ell_1)\otimes_i C_1(\max))\otimes_j(B_\alpha(\ell_2)\otimes_iD_\beta(\ell_3))\\
&(A_1(\max)\otimes_i C_2(\ell_1))\otimes_j(B_\alpha(\ell_2)\otimes_iD_\beta(\ell_3))\\
&(A_\alpha(\ell_1)\otimes_iC_\beta(\ell_2))\otimes_j(B_2(\ell_3)\otimes_i D_1(\max))\\
&(A_\alpha(\ell_1)\otimes_iC_\beta(\ell_2))\otimes_j(B_1(\max)\otimes_i D_2(\ell_3))
\end{aligned}
\end{equation}
Here $\max$ is the maximal possible value of the current parameter.
These elements are not necessarily non-vanishing, but for some $\alpha,\beta,\ell_s$ they (and only they) may not vanish.

Now some of these elements will necessarily vanish by the third map.

{\bf Regular elements which may not vanish by the third map:}
\begin{equation}\label{keyvanishing3}
(A_\alpha(\ell_1)\otimes_iC_\alpha(\ell_2))\otimes_j(B_\alpha(\ell_3)\otimes_iD_\alpha(\ell_4))
\end{equation}

{\bf Singular elements which may not vanish by the third map:}

\begin{equation}\label{keyvanishing4}
\begin{aligned}
\ &(A_1(\max)\otimes_iC_1(\max))\otimes_j(B_2(\ell_1)\otimes_iD_2(\ell_2))\\
&(A_2(\ell_1)\otimes_i C_2(\ell_2))\otimes_j(B_1(\max)\otimes_i D_1(\max))
\end{aligned}
\end{equation}

{\bf Conclusion:}

We see from \eqref{keyvanishing3} that the only regular non-vanishing elements are exactly those for which all indices $\alpha,\beta,\gamma,\delta$ are equal (exactly as is stated by Lemma).

Now consider the non-vanishing singular elements.

We distinguish {\it simple} singular elements (those for which only one of the two last maps of the right-hand path of \eqref{compexpl} are in singular position),
and singular elements {\it of the second order} (which are singular for the both maps).

{\it Simple singular elements:}

Suppose we have an element \eqref{keyvanishingelement} for which the second map is regular and the third is singular. Then clearly it is of the form \eqref{keyvanishing4}, and only them (totally 2 elements). They are corresponded to the second and to the fifth lines in \eqref{keyvanishingfinal2}.

Suppose we have an element \eqref{keyvanishingelement} for which the second map is singular and the third is regular. This case is corresponded to the elements \eqref{keyvanishing2} with either $\alpha=\beta=2$, or one of indices $\alpha,\beta$ is equal to 1 and the corresponding $\ell_s=\max$, and the remaining index is 2. In the both cases the third map is regular. Indeed, the elements $X_1(\max)$ $(X=A,B,C,D)$ with the lower index 1 will be mapped to $X_2(\ell)$ (for $\ell=1$, but it is not essential). Then we get an element of the form
\eqref{keyvanishing3} with $\alpha=2$ which is regular for the third map. The first alternative gives the whole group \eqref{keyvanishingfinal1}, while the second alternative gives the remaining 4 elements (all lines except the second and the fifth) in \eqref{keyvanishingfinal2}.

{\it Singular elements of the second order:}

To get singular elements for the third maps \eqref{keyvanishing4} from the singular elements for the second map \eqref{keyvanishing2} we necessarily should have in \eqref{keyvanishing2} $\alpha=\beta=1$, and the corresponding two $\ell_s=\max$. This gives the whole third group \eqref{keyvanishingfinal3}.

The Lemma is proven for the right-hand path.

The analysis of the left-hand path is fairly analogous, and we leave it to the reader.

\endproof
Now we can prove
\begin{lemma}
The diagram \eqref{compexpl} is commutative.
\end{lemma}
\proof{}
We know from Lemma \ref{spisok} the list of all elements on which the left-hand path and the right-hand path of \eqref{compexpl} do not vanish, and we know that this list is precisely the same for the both pathes. Now the proof is a straightforward check that the both pathes coincide on any of these elements.

For the regular elements (which are the same for the both pathes) the statement is clear (it follows from the diagrams \eqref{intexpl} and \eqref{extexpl} in the underlying $n$-fold monoidal category $\mathcal{C}$). For the singular elements \eqref{keyvanishingfinal1}-\eqref{keyvanishingfinal3} it is a direct and a routine check; we leave it to the reader avoiding to make longer already a very long and technical proof of Key-Lemma \ref{extmon}.

\endproof

The Key-Lemma \ref{extmon} is proven.
\endproof

\comment
\subsection{Modules over $n$-fold monoidal categories}
Here we develop some mini-theory which will be applied in Section 4
for an elegant proof that tetramodules form a 2-fold monoidal
category. More generally, it can be applied to a proof that the
tetramodules over $n$-fold monoidal algebras (see Example
\ref{monalg}) form an $(n+1)$-monoidal category.

\subsubsection{The case of usual categories}
Let $\mathcal{C}$ be a $k$-linear category. We say that a {\it
module} over it is the set of $k$-vector spaces
$\{M_{X_1,X_2}\}_{X_1,X_2\in \mathrm{Ob}(\mathcal{C})}$ such that if
we set $\Mor_M(X_1,X_2)=\Mor(X_1,X_2)\oplus \epsilon M_{X_1,X_2}$,
where $\epsilon^2=0$, we get a category $\mathcal{C}_M$, such that
the natural inclusion $i_M\colon \mathcal{C}\to\mathcal{C}_M$ is a
map of categories.
\begin{example}\label{example0}
Let our category have a single object, $a$, such that $\Mor(a,a)=A$.
Then $A$ is an associative algebra, and a module over this category
is the same that a bimodule over $A$.
\end{example}

In general, $\Mor_M(X_1,X_2)$ should be a right
$\Mor_\mathcal{C}(X_1,X_1)$-module and a left
$\Mor_\mathcal{C}(X_2,X_2)$-module; these two module structures
commute. Actually, a structure of a module over the category
$\mathcal{C}$ is the same that the set of
$\Mor_\mathcal{C}(X_2,X_2)$-$\Mor_\mathcal{C}(X_1,X_1)$-bimodules.

\begin{example}\label{example1}
Let $A$ be an associative algebra, and let $\mathcal{C}$ be the
category of left $A$-modules. Let $M$ be an $A$-bimodule. Define
\begin{equation}\label{modalg}
\Mor_M(X_1,X_2)=\Mor_\mathcal{C}(M\otimes_A X_1, M\otimes _A X_2)
\end{equation}
Then this set of morphism gives a module over the category of left
$A$-modules, depending on an $A$-bimodule $M$.
\end{example}

For any two modules over a $k$-linear category $\mathcal{C}$ $M,N$
one can define their {\it tensor product} $M\otimes N$. By
definition,

\begin{equation}\label{eqmod1}
\begin{aligned}
\ &(M\otimes N)_{X_1,X_2}=M_{X_1,X_2}\otimes_k N_{X_1,X_2}/\\&
\bigl\{m.(f_1f_2)\otimes n-m\otimes (f_1f_2).n=0\text{ for any
$f_1\in\Mor_\mathcal{C}(X_1,X_2)$ and
$f_2\in\Mor_\mathcal{C}(X_2,X_1)$}\bigl\}
\end{aligned}
\end{equation}

\begin{example}
Consider the situation of the Example \ref{example0}. Then the
tensor product of modules over categories, corresponding to
$A$-bimodules $M_1$ and $M_2$, is the $A$-bimodule $M_1\otimes_A
M_2$.
\end{example}

\begin{example}
Let $\mathcal{C}$ be the category of free left $A$-modules, $A$ an
associative algebra. Consider $A$-bimodules $M$ and $N$ which a free
as left and free as right modules. Then each of them defines a
module over $\mathcal{C}$, by \eqref{modalg}. Their tensor product
comes by the same formula from $M\otimes_A N$.
\end{example}

\begin{lemma}
Let $\mathcal{C}$ be a $k$-linear category. Then the modules over
$\mathcal{C}$ form a monoidal category, with the monoidal structure
defined by \eqref{eqmod1}.
\end{lemma}
\endproof
\endcomment

\section{The category of tetramodules}
\subsection{}
Recall that an {\it associative bialgebra} is a vector space $A$ over a field $k$
equipped with two operations, the product $*:A^{\otimes 2}\to A$ and
the coproduct $\Delta\colon A\to A^{\otimes 2}$, which obey the
axioms 1.-4. below:
\begin{itemize}
\item[1.] Associativity: $a*(b*c)=(a*b)*c$;
\item[2.] Coassociativity: $(\Delta\otimes
id)\Delta(a)=(id\otimes\Delta)\Delta(a)$;
\item[3.] Compatibility: $\Delta(a*b)=\Delta(a)*\Delta(b)$.
\end{itemize}
We use the classical notation
$$
\Delta(a)=\Delta^1(a)\otimes\Delta^2(a)
$$
which is just a simplified form of the equation
$$
\Delta(a)=\sum_i\Delta^1_i(a)\otimes\Delta^2_i(a)
$$
We always assume that our bialgebras have a unit and a counit.
A unit is a map $i\colon k\to A$ and the counit is a map $\varepsilon\colon A\to k$. We always assume
\begin{itemize}
\item[4.] $i(k_1\cdot k_2)=i(k_1)*i(k_2)$, $\varepsilon(a*b)=\varepsilon(a)\cdot \varepsilon(b)$.
\end{itemize}
We also denote the product $*$ by $m$.

A {\it Hopf algebra} is a bialgebra with antipode. An antipode is a $k$-linear map $S\colon A\to A$ which obeys
\begin{itemize}
\item[5.] $m(1\otimes S)\Delta(a)=m(S\otimes 1)\Delta(a)=i(\varepsilon(a))$
\end{itemize}
We need the existence of antipode in only one, but a crucial place, in Section 4.3. All results of Sections 4.1 and 4.2 do not need the existence of antipode.

As we already mentioned here, an ``operadic'' definition of a
bimodule over an associative algebra $A$ is a $k$-vector space $M$
such that $A\oplus\epsilon M$ is again associative algebra, where
$\epsilon^2=0$, and the restriction of the algebra structure to $A$
coincides with the initial one. Such $M$ is the same that an
$A$-bimodule. We give an analogous definition in the case when $A$
is an associative bialgebra.

\begin{defn}
Let $A$ be an associative bialgebra. A Bernstein-Khovanova
tetramodule $M$ over $A$ is a vector space such that $A\oplus
\epsilon M$ is an associative bialgebra when $\epsilon^2=0$ and the
restriction of the bialgebra structure to $A$ is the initial one.
The category of tetramodules over a bialgebra $A$ is denoted
$\Tetra(A)$.
\end{defn}

More precisely, one has maps $m_\ell\colon A\otimes M\to M$,
$m_r\colon M\otimes A\to M$ (which make $M$ an $A$-bimodule), and
maps $\Delta_\ell\colon M\to A\otimes M$ and $\Delta_r\colon M\to
M\otimes A$ (which make $M$ an $A$-bicomodule), with some
compatibility between these 4 maps. The compatibility written
explicitly is the following:
\begin{equation}\label{eq4.10}
\Delta_\ell(a*m)=(\Delta^1(a)*\Delta_\ell^1(m))\otimes
(\Delta^2(a)*\Delta^2_\ell(m))\subset A\otimes_k M
\end{equation}

\begin{equation}\label{eq4.11}
\Delta_\ell(m*a)=(\Delta_\ell^1(m)*\Delta^1(a))\otimes
(\Delta_\ell^2(m)*\Delta^2(a))\subset A\otimes_k M
\end{equation}

\begin{equation}\label{eq4.12}
\Delta_r(a*m)=(\Delta^1(a)*\Delta_r^1(m))\otimes
(\Delta^2(a)*\Delta^2_r(m))\subset  M\otimes_k A
\end{equation}

\begin{equation}\label{eq4.13}
\Delta_r(m*a)=(\Delta_r^1(m)*\Delta^1(a))\otimes
(\Delta^2_r(m)*\Delta^2(a))\subset M\otimes_k A
\end{equation}
Here we use the natural notation like
$\Delta_\ell(m)=\Delta^1_\ell(m)\otimes\Delta^2_\ell(m)$ with
$\Delta_\ell^1(m)\in A$, $\Delta^2_\ell(m)\in M$, etc. As well, we
use the sign $*$ for the both product in $A$ and the module products
$m_\ell$ and $m_r$.

The main example of a tetramodule over $A$ is a itself; it is called
the tautological tetramodule.

When $A$ is finite-dimensional over $k$, a tetramodule is the same
that a left module over some associative algebra $H(A)$. This
algebra $H(A)$ is, as a vector space, the tensor product
$H(A)=A\otimes_kA\otimes_kA^*\otimes A^*$, and the commutation
relation are such that the equations \eqref{eq4.10}-\eqref{eq4.13}
above agree. (This algebra $H(A)$ is the ``double Heisenberg
double'' of the bialgebra $A$).

In particular, if $A$ is finite-dimensional over $k$, the category
$\Tetra(A)$ has enough projectives and enough injectives objects.
For general $A$, R.Taillefer proved [Tai2] that the category
$\Tetra(A)$ has enough injectives.

The main relation with the deformation theory is the following
theorem, proven by R.Taillefer [Tai1,2]:

\begin{theorem}
Let $H^\mb_\GS(A,A)$ denote the Gerstenhaber-Schack cohomology of an
associative bialgebra $A$. Then one has:
\begin{equation}\label{eq4.14}
H^\mb_\GS(A,A)=\Ext^\mb_{\Tetra(A)}(A,A)
\end{equation}
\end{theorem}
\endproof

The Gerstenhaber-Schack cohomology $H^\mb_\GS(A,A)$ is known to
control the infinitesimal deformations of the bialgebra $A$ [GS]. We overview the Gerstenhaber-Schack cohomology
and prove Theorem 4.2 in the Appendix Section 6.

\begin{remark}
To control the global deformations of a bialgebra $A$ through the
usual deformation theory (the Maurer-Cartan equation etc.) one needs
to have an appropriate $L_\infty$ structure on the
Gerstenhaber-Schack complex. (In this case it will be an $L_\infty$
structure with the components at least up to 4th, because the r.h.s.
of the structure equation of bialgebra,
$\Delta(a*b)=\Delta(a)*\Delta(b)$, is of 4th degree in the
operations. This structure is not known yet. What we construct in
Section 4.3 is a Lie bracket which presumably is induced on the
cohomology from this structure.
\end{remark}

\begin{example}
Consider the case when $A=S(V)$ is a free (co)commutative bialgebra,
for simplicity suppose $V$ is finite-dimensional over $k$. We prove in Section 6 that
\begin{equation}\label{eq4.15}
H^k_\GS(A,A)=\oplus_{i+j=k}\Lambda^iV\otimes_k\Lambda^jV^*
\end{equation}
The wedge-product defines on the r.h.s. a super-commutative product
of degree 0, and the contraction of $V$ and $V^*$ defines a Lie
bracket of degree -2. Clearly the bracket is Poisson and even, it
does not cause an additional sign. The product and the bracket obey
the even Leibniz compatibility. Altogether it defines {\it a
3-algebra structure} on $H^\mb_\GS(A,A)$ for $A=S(V)$. In this paper
we construct such a structure for any bialgebra $A$.
\end{example}

\subsection{The structure of a 2-fold monoidal category on
$\Tetra(A)$}

\subsubsection{Two ``external'' tensor products}
 Define firstly for a pair of $A$-tetramodules $M_1,M_2$ two their
``external'' tensor products $M_1\boxtimes_1M_2$ and
$M_1\boxtimes_2M_2$ which are again $A$-tetramodules. In the both
cases the underlying vector space is $M_1\otimes_kM_2$.

{The case of $M_1\boxtimes_1M_2$}:
\begin{itemize}
\item[1.] $m_\ell(a\otimes m_1\boxtimes m_2)=(am_1)\boxtimes m_2$,
\item[2.] $m_r(m_1\boxtimes m_2\otimes a)=m_1\boxtimes (m_2a)$,
\item[3.] $\Delta_\ell(m_1\boxtimes
m_2)=(\Delta_\ell^1(m_1)*\Delta_\ell^1(m_2))\otimes
(\Delta^2_\ell(m_1)\boxtimes\Delta^2_\ell(m_2))$,
\item[4.] $\Delta_r(m_1\boxtimes m_2)=(\Delta^1_r(m_1)\boxtimes
\Delta^1_r(m_2))\otimes (\Delta^2_r(m_1)*\Delta^2_r(m_2))$.
\end{itemize}

{The case of $M_1\boxtimes_2 M_2$}:
\begin{itemize}
\item[1.] $m_\ell(a\otimes m_1\boxtimes
m_2)=(\Delta^1(a)m_1)\boxtimes (\Delta^2(a)m_2)$,
\item[2.] $m_r(m_1\boxtimes m_2\otimes
a)=(m_1\Delta^1(a))\boxtimes(m_2\Delta^2(a))$,
\item[3.] $\Delta_\ell(m_1\boxtimes
m_2)=\Delta^1_\ell(m_1)\otimes(\Delta_\ell^2(m_1)\boxtimes m_2)$,
\item[4.] $\Delta_r(m_1\boxtimes
m_2)=(m_1\boxtimes\Delta_r^1(m_2))\otimes\Delta_r^2(m_2)$.
\end{itemize}

For the both definitions we do not use the whole tetramodule
structures on $M_1,M_2$. In particular, for the first definition we
never use the right multiplication $m_r$ for $M_1$ and the left
multiplication $m_\ell$ for $M_2$. As well, for the second
definition we do not use $\Delta_r$ for $M_1$ and $\Delta_\ell$ for
$M_2$.

This gives us some additional possibilities, which we use to define
the ``internal'' tensor products $M_1\otimes_1 M_2$ and
$M_1\otimes_2 M_2$. For the both cases the tautological tetramodule
is a unit object.

\subsubsection{Two ``internal'' tensor products}
\begin{defn}
Let $M_1,M_2$ be two tetramodules over a bialgebra $A$. Their first
tensor product $M_1\otimes_1 M_1$ is defined as the
quotient-tetramodule
\begin{equation}\label{mon1}
M_1\otimes_1M_2=M_1\boxtimes_1M_2/((m_1a)\boxtimes_1m_2-m_1\boxtimes_1(am_2))
\end{equation}
One easily checks that this definition is correct. Analogously, the
second tensor product $M_1\otimes_2M_2$ is defined as a
sub-tetramodule
\begin{equation}\label{mon2}
M_1\otimes_2M_2=\{\sum_im_{1i}\boxtimes_2 m_{2i}\subset
M_1\boxtimes_2 M_2|\sum_i \Delta_r(m_{1i})\otimes_k
m_{2i}=\sum_im_{1i}\otimes_k\Delta_\ell(m_{2i})\}
\end{equation}
Again, one easily checks that this definition is correct.
\end{defn}
\begin{lemma}\label{units}
Suppose that the bialgebra $A$ has a unit and a counit. Then the
tautological tetramodule $A$ is the unit object for the both
monoidal structures.\end{lemma} \proof{} Let $M$ be a tetramodule.
One can check that the following maps
\begin{equation}\label{a1}
\begin{aligned}
\ & m_\ell\colon A\otimes_1M\to M\\
&m_r\colon M\otimes_1A\to M
\end{aligned}
\end{equation}
and
\begin{equation}\label{a2}
\begin{aligned}
\ & \Delta_\ell\colon M\to A\otimes_2M\\
&\Delta_r\colon M\to M\otimes_2A
\end{aligned}
\end{equation}
are morphisms of tetramodules. If $A$ has a unit, the first two maps
are isomorphisms, while if $A$ has a counit, the second ones two
are.
\endproof

\subsubsection{The 2-fold monoidal structure}
We construct the map $\eta_{M,N,P,Q}\colon (M\otimes_2
N)\otimes_1(P\otimes_2 Q)\to (M\otimes_1 P)\otimes_2(N\otimes_1 Q)$
in several steps.

{\bf The first step} is to check that the map $\phi_0\colon
(M\boxtimes_2 N)\boxtimes_1(P\boxtimes_2 Q)\to (M\boxtimes_1
P)\boxtimes_2 (N\boxtimes_1 Q)$, $\phi_0(m\otimes_k n\otimes_k
p\otimes_k q)=m\otimes_k p\otimes_k n\otimes_k q$, is a map of
tetramodules. We have:
\begin{equation}\label{eqformula1}
\begin{aligned}
\ &a*\bigl((m\boxtimes_2 n)\boxtimes_1 (p\boxtimes_2
q)\bigr)=\\
&\bigl(a*(m\boxtimes_2 n)\bigr)\boxtimes_1 (p\boxtimes_2
q)=\\
&(\Delta^1(a)*m)\otimes(\Delta^2(a)*n)\otimes p\otimes q
\end{aligned}
\end{equation}
and
\begin{equation}\label{eqformula2}
\begin{aligned}
\ &a*\bigl((m\boxtimes_1 p)\boxtimes_2(n\boxtimes_1 q)\bigr)=\\
&\bigl(\Delta^1(a)*(m\boxtimes_1
p)\bigr)\boxtimes_2\bigl(\Delta^2(a)*(n\boxtimes_1 q)\bigr)=\\
&(\Delta^1(a)*m)\otimes p\otimes (\Delta^2(a)*n)\otimes q
\end{aligned}
\end{equation}
We see that
\begin{equation}
\phi_0(\text{r.h.s. of \eqref{eqformula1}})=\text{r.h.s. of
\eqref{eqformula2}}
\end{equation}
That is, $\phi_0$ is a map of left modules; analogously it is a map
of right modules.

Now prove that $\phi_0$ is a map of left comodules. We have:
\begin{equation}\label{eqformula3}
\begin{aligned}
\ &\Delta_\ell\bigl((m\boxtimes_2 n)\boxtimes_1 (p\boxtimes_2
q)\bigr)=\\
&\bigl(\Delta_\ell^1(m\boxtimes_2 n)*\Delta_\ell^1(p\boxtimes_2
q)\bigr)\otimes_k\bigl(\Delta_\ell^2(m\boxtimes_2
n)*\Delta_\ell^2(p\boxtimes_2 q)\bigr)=\\
&\bigl(\Delta_\ell^1(m)*\Delta_\ell^1(p)\bigr)\otimes_k\bigl(\Delta_\ell^2(m)\otimes_k
n\otimes_k\Delta_\ell^2(p)\otimes_k q\bigr)
\end{aligned}
\end{equation}
and
\begin{equation}\label{eqformula4}
\begin{aligned}
\ &\Delta_\ell\bigl((m\boxtimes_1 p)\boxtimes_2(n\boxtimes_1
q)\bigr)=\\
&\Delta_\ell^1(m\boxtimes_1p)\otimes_k\bigl(\Delta_\ell^2(m\boxtimes_1
p)\boxtimes_2(n\boxtimes_1 q)\bigr)=\\
&\bigl(\Delta_\ell^1(m)*\Delta_\ell^1(p)\bigr)\otimes_k\bigl(\Delta_\ell^2(m)\otimes_k
\Delta_\ell^2(p)\otimes_k n\otimes_k q\bigr)
\end{aligned}
\end{equation}
We see that
\begin{equation}
\Delta_\ell\circ \phi_0=id\otimes_k(\phi_0\circ\Delta_\ell)
\end{equation}
that is, $\phi_0$ is a map of left comodules. It is proven
analogously that $\phi_0$ is a map of right comodules.

\vspace{2mm}

{\bf At the second step} we consider the natural projections of
tetramodules $p_{M,P}\colon M\boxtimes_1 P\to M\otimes_1 P$ and
$p_{N,Q}\colon N\boxtimes_1Q\to N\otimes_1 Q$. We consider the
composition
\begin{equation}
\phi_1=(p_{M,P}\boxtimes_2 p_{N,Q})\circ\phi_0\colon (M\boxtimes_2
N)\boxtimes_1(P\boxtimes_2 Q)\to (M\otimes_1
P)\boxtimes_2(N\otimes_1 Q)
\end{equation}
We want to check that the map $\phi_1$ defines naturally a map
\begin{equation}
\phi_2=\overline{\phi}_1\colon (M\boxtimes_2
N)\otimes_1(P\boxtimes_2 Q)\to (M\otimes_1 P)\boxtimes_2(N\otimes_1
Q)
\end{equation}
that is, the elements of the form
\begin{equation}\label{eqformula10}
\bigl((m\boxtimes_2 n)*a\bigr)\otimes_k (p\boxtimes_2
q)-(m\boxtimes_2 n)\otimes_k \bigl(a*(p\boxtimes_2 q)\bigr)
\end{equation}
are mapped to 0 by $\phi_1$.

Indeed,
\begin{equation}
\eqref{eqformula10}=(m*\Delta^1a)\otimes_k(n*\Delta^2a)\otimes_k(p\otimes_k
q)-(m\otimes_k n)\otimes_k (\Delta^1a*p)\otimes_k (\Delta^2a*q)
\end{equation}
which, after the permutation $\phi_0$ of the second and the third
factors, is mapped to 0 in $(M\otimes_1 P)\boxtimes_2(N\otimes_1
Q)$. Therefore, the map $\phi_2$ is well-defined.

\vspace{2mm}

{\bf At the third step} we restrict the map $\phi_2$ to $(M\otimes_2
N)\otimes_1(P\otimes_2 Q)\subset (M\boxtimes_2 N)\otimes_1
(P\boxtimes_2 Q)$, and we need to check that the image of this
restricted map belongs to $(M\otimes_1 P)\otimes_2(N\otimes_1
Q)\subset (M\otimes_1 P)\boxtimes_2(N\otimes_1 Q)$.

Suppose $m\otimes_k n\in M\otimes_2N$ and $p\otimes_k q\in
P\otimes_2 Q$ (we assume the summation over several such monomials,
but for simplicity we skip this summation). Then
\begin{equation}\label{eqformula11}
\Delta_r^1(m)\otimes_k\Delta_r^2(m)\otimes_k
n=m\otimes_k\Delta_\ell^1(n)\otimes_k\Delta_\ell^2(n)
\end{equation}
with the middle factors in $A$, and analogously
\begin{equation}\label{eqformula12}
\Delta_r^1(p)\otimes_k\Delta_r^2(p)\otimes_kq=p\otimes_k\Delta_\ell^1(q)\otimes_k\Delta_\ell^2(q)
\end{equation}
again, with the middle factors in $A$.

One needs to prove that \eqref{eqformula11} and \eqref{eqformula12}
together imply that
\begin{equation}\label{eqformula13}
(m\boxtimes_1 p)\boxtimes_2(n\boxtimes_1 q)\in (M\otimes_1
P)\otimes_2 (N\otimes_1 Q)\subset(M\otimes_1 P)\boxtimes_2
(N\otimes_1 Q)
\end{equation}
that is,
\begin{equation}\label{eqformula14}
\Delta_r^1(m)\otimes_k
\Delta_r^2(p)\otimes_k(\Delta_r^2(m)*\Delta_r^2(p))\otimes_k
n\otimes_kq=m\otimes_k p\otimes_k
(\Delta_\ell^1(n)*\Delta_\ell^1(q))\otimes_k
\Delta_\ell^2(n)\otimes_k\Delta_\ell^2(q)
\end{equation}

To get \eqref{eqformula14} from \eqref{eqformula12} and
\eqref{eqformula13} we permute \eqref{eqformula12} such that the
factors in $A$ are the most right, permute \eqref{eqformula12} such
that the factors in $A$ are the most left, then take the equation
$(\text{l.h.s. of }\eqref{eqformula12})\otimes_k(\text{l.h.s. of
}\eqref{eqformula13})=(\text{r.h.s. of
}\eqref{eqformula12})\otimes_k(\text{r.h.s. of
}\eqref{eqformula13})$ (after they are permuted). Then for the two
middle factor (in $A$) we apply the product $*\colon A\otimes_k A\to
A$, and then permute again.

\vspace{3mm}

{\bf The map $\mathbf{\eta_{M,N,P,Q}}$ is constructed.} \vspace{2mm}

\begin{theorem}
The maps $\eta_{M,N,P,Q}$ constructed above, together with the two
tensor products $\otimes_1$ and $\otimes_2$, define a 2-fold
monoidal structure on the category $\Tetra(A)$ of tetramodules over
a bialgebra $A$.
\end{theorem}
\proof{} First of all, the two tensor products $\boxtimes_1$ and
$\boxtimes_2$ with $\wtilde{\eta}_{M,N,P,Q}=\phi_0$, clearly forms a
2-fold monoidal structure on the category $\Tetra(A)$ (but this
2-fold monoidal structure does not admit unit elements). In
particular, the diagrams \eqref{intexpl} and \eqref{extexpl} are
commutative for $\wtilde{\eta}_{M,N,P,Q}$ (because $\phi_0$ is just
the permutation which switches the second and the third factors).
Now the same diagrams for the actual structure $\eta_{M,N,P,Q}$ are
obtained from these simple ones by passing to subquotients.
Therefore, they are commutative as well.
\endproof

\subsection{The condition (*) for $\Tetra(A)$}
Recall the meaning of the condition (*): we want to restrict ourselves with a homotopy equivalent to $\EXT^k(A,A)$ subcategories on which the monoidal structures are {\it exact}. This condition is essential in the Schwede's construction described in Section 2, and in its generalization described in the sequel Section 4.4.

We have two monoidal structures on $\Tetra(A)$, namely $\otimes_1$ and $\otimes_2$. We leave to the reader the following simple statement:
\begin{lemma}
The monoidal bifunctor $\otimes_1$ is right exact, and the monoidal bifunctor $\otimes_2$ is left exact.
\end{lemma}
\endproof

The author does not know, in this setting, how to construct a homotopically equivalent subcategory on which the both monoidal structures are exact, if we work with general bialgebras.
It seems that the construction of [Sch], Lemma 2.1 can not be adopted to our situation.

But when we suppose that our bialgebra is a Hopf algebra, that is, it obeys an antipode, the situation is much better because of the following classical result:
\begin{lemma}
Let $A$ be a Hopf algebra (that is, a bialgebra with an antipode), and let $M$ be a $k$-vector space with structures of left $A$-module
$m_\ell\colon A\otimes M\to M$, and of left $A$-comodule $\Delta_\ell\colon M\to A\otimes M$, which are compatible as
\begin{equation}
\Delta_\ell(a*m)=(\Delta^1(a)*\Delta_\ell^1(m))\otimes (\Delta^2(a)*\Delta_\ell^2(m))
\end{equation}
(this is just a $\frac14$ of the structure of a tetramodule). Then $M$ is necessarily free as $A$-module and cofree as $A$-comodule.
\end{lemma}
\proof{}
See [Sw], Section 4.1 for a proof of the first claim; the second claim is dual to the first one.
\endproof
We see from the Lemma above that the both monoidal structures $\otimes_1$ and $\otimes_2$ are exact on the nose, so there is no necessity to restrict by a smaller category.

For this we should include the existence of antipode condition to all statements which use the condition (*). That is, in all such cases we should work not with general bialgebras, but with Hopf algebras.

\subsection{A generalization of the Schwede's construction}
Let $A$ be a Hopf algebra. Here we construct a Lie bracket of degree -2 on the graded space $\Ext^\mb_{\Tetra(A)}(A,A)$ generalizing the Schwede's construction of the Lie bracket of degree -1 on the Hochschild cohomology, see Section 2.2.

Let $M,N$ be two extensions of tetramodules, $M\in \EXT^m(A,A)$, $N\in\EXT^n(A,A)$. Consider the following ``big octahedron'' diagram. We prove just below in Lemma \label{commtriangles} that the triangle 2-faces in \eqref{maindiagram} are commutative. That is, the diagram defines an element in $\pi_2(\EXT^{m+n}(A,A))$. By Retakh's theory (see Section 1) the latter is isomorphic to $\Ext^{m+n-2}(A,A)$. This is the construction of the bracket. A priori it is not clear that this operation satisfies the Jacobi identity, it follows from much deeper results of Section 5.

\begin{equation}\label{maindiagram}
\xymatrix{ &&M\sharp N\\
\\
&&&N\otimes_{2,\tau}M\ar[luu]\ar[ddddl]\\
M\otimes_{1,\tau}N\ar[urrr]\ar[dr]\ar[rruuu]\ar[rrddd]&&&&N\otimes_{1,\tau}M\ar[ul]\ar[dlll]\ar[lluuu]\ar[llddd]\\
&M\otimes_{2,\tau}N\ar[ruuuu]\ar[rdd]\\
\\
&&N\sharp M }
\end{equation}

\begin{lemma}\label{commtriangles}
All triangle 2-faces in \eqref{maindiagram} are commutative diagrams.
\end{lemma}
\proof{}
We prove a more general statement. Given an $n$-fold monoidal category, let $1\le i<j<k\le n$.
Consider the maps $\varphi_{ij}\colon M\otimes_iN\to M\otimes_jN$, $\varphi_{jk}\colon M\otimes_jN\to M\otimes_kN$, and $\varphi_{ik}\colon M\otimes_iN\to M\otimes_k N$, defined in \eqref{12left}. Recall that $\varphi_{ij}=\eta^{ij}_{M, A,A,N}$ where $A$ is the common unit object.
We prove that
\begin{equation}\label{commtri}
\varphi_{ik}=\varphi_{jk}\circ \varphi_{ij}
\end{equation}
This is an application of the commutative diagram \eqref{compexpl} (this diagram is indeed commutative for $n=3$ and $\otimes_3=\sharp$ by the Key-Lemma \ref{extmon}). Consider the diagram \eqref{compexpl} for $A_1=M$, $A_2=B_1=B_2=C_1=C_2=D_1=A$, $D_2=N$. Then the right-hand path of \eqref{compexpl} gives $\varphi_{jk}\circ \varphi_{ij}$ while the left-hand path gives $\varphi_{ik}$, and \eqref{commtri} follows.
This proves the commutativity of one of the 8 triangle faces in \eqref{maindiagram}. The commutativity of the remaining 7 triangles is proven analogously, using $\theta_{ij}\colon M\otimes_i N\to N\otimes_j M$ defined in \eqref{12right}, and different substitutions to the diagram \eqref{compexpl}.
\endproof

\begin{remark}
This lemma, and the commutativity of the triangle faces in \eqref{maindiagram}, is a nice example on (a rather complicated) diagram \eqref{compexpl}, and a motivation for a rather technical proof of the Key-Lemma \ref{extmon}.
\end{remark}

\begin{remark}
It follows from the existence of the big octahedron that the ``bracket'' of degree -1 defined as in Section 2.2 either from $\otimes_1$ or $\otimes_2$ separately, is equal to 0. In particular, the Lie bracket defined in [Ta2], Section 5, is identically 0. The idea is that in the presence of two monoidal structures compatible as in 2-fold monoidal category, the Schwede's loop \eqref{eq2.5} is ``divided'' by 4 commutative triangles, and therefore is contractible in the nerve. Morally, to get non-trivial operations, we should take into account all possible monoidal structures. 
\end{remark}

\section{Passage to Spectra}

\subsection{From ``spectra'' of categories to spectra of topological spaces}
\subsubsection{Appearance of a proplem}
There is an operad of categories $\{\mathcal{M}_{n+1}(d)\}$ acting on each $(n+1)$-fold monoidal category and, in particular, on $\bigsqcup_{k\ge 0}
\EXT^k_{\mathcal{C}}(A,A)$, where $\mathcal{C}$ is an $n$-fold abelian monoidal category satisfying the condition (*), and $A$ is a common unit object in it (see Section 3). Moreover, there are ``spectrum structure maps'' $G_k\colon \EXT^k(A,A)\to\Omega_R\EXT^{k+1}(A,A)$ (see Section 1.3). Our goal is to pass, by the classifying space functor, to an $\Omega$-spectrum $X$ of topological spaces with $X_k=|\EXT^k(A,A)|$ such that the operad of spectra $\Sigma^\infty|\mathcal{M}_{n+1}(d)|$ acts on this spectrum (after suitable reducing of the both spectra to $\Sigma$-spectra).

This compatibility of the operad action with the structure maps of spectra is more natural (and, seemingly, only possible) to prove on the level of categories. Here the problem we meet is the following. It is natural to consider the spectra of {\it based} topological spaces. Therefore, the corresponding spectrum of categories should be also considered as based. But our operad is {\it not} compatible with any sense of the based objects, in the sense that the $n+1$ monoidal structures do not give a based object if one of two its arguments is a based object.

More precisely, we would like to prove the commutativity of the following diagram of functors:
\begin{equation}\label{eq5.1}
\xymatrix{
\EXT^k\times\EXT^\ell\ar[r]^{\otimes_i}\ar[d]_{G_k\otimes id}&\EXT^{k+\ell}\ar[d]^{G_k}\\
(\Omega_R\EXT^{k+1})\times \EXT^\ell\ar[r]^{\otimes_i}&\Omega_R\EXT^{k+\ell+1}}
\end{equation}
However, this diagram does not make sense, because the lower horizontal arrow is ill-defined. Indeed, if $*\rightarrow X\leftarrow *$ is and element in $\Omega_R$, we would like to induce from any monoidal structure the component-wise product with $Y$, which should be $*\otimes Y\rightarrow X\otimes Y\leftarrow *\otimes Y$. The point is that this element is not an object of $\Omega_R$ because $*\otimes Y\ne *$.

The origin of this problem is that the $n+1$ monoidal structures $\otimes_i\colon \EXT^a(A,A)\times \EXT^b(A,A)\to \EXT^{a+b}(A,A)$ (defined for an $n$-fold monoidal abelian category $\mathcal{C}$), {\it do not respect the based objects}. Consequently, passing to the nerves we get a map
$\otimes_i\colon B(\EXT^a(A,A))\times B(\EXT^b(A,A))\to B(\EXT^{a+b}(A,A))$, and {\it not} a map
$\otimes_i\colon B(\EXT^a(A,A))\wedge B(\EXT^b(A,A))\to B(\EXT^{a+b}(A,A))$ (where $\wedge$ is as usual the direct product in the category of based topological spaces). This circumstance makes impossible to use directly the smash-products of spectra, which we wish to use, in order to pass to the operad action on homotopy groups of the spectrum $\{B(\EXT^k(A,A))\}$.

\comment
\subsubsection{Our strategy to solve it}
Let $\mathcal{C}$ be a category and let $\mathcal{D}\subset\mathcal{C}$ be a subcategory. We define two new categories,
$\mathbf{C}$ and $\mathbf{D}$.

A {\it path} of length $n$ in the category $\mathcal{C}$ is a chain of maps
\begin{equation}\label{chains}
X_0\rightarrow X_1\leftarrow X_2 \rightarrow X_3\leftarrow\dots \leftrightarrow X_n
\end{equation}
where all $X_i$ are objects of $\mathcal{C}$. A morphism of two pathes is defined as the corresponding commutative diagrams.
Denote this category by $\mathcal{P}ath_n(\mathcal{C})$.
There are canonical functors $\mathcal{P}ath_n(\mathcal{C})\to \mathcal{P}ath_m(\mathcal{C})$ for $n\le m$ (one set $Y_i=X_i$ for $1\le i\le n$, $Y_i=X_n$ for $n+1\le i\le m$, and the new maps are identities).

Define the category $\mathbf{C}$ as the inductive limit
\begin{equation}
\mathbf{C}=\underrightarrow{\lim}\mathcal{P}ath_n(\mathcal{C})
\end{equation}

Analogously, define the category $\mathcal{P}ath_n(\mathcal{C}_\mathcal{D})$ whose objects are chains \eqref{chains} with all $X_i\in\mathcal{C}$ but $X_0\in\mathcal{D}$.

Again, set
\begin{equation}
\mathbf{D}=\underrightarrow{\lim}\mathcal{P}ath_n(\mathcal{D}_\mathcal{D})
\end{equation}

We have the following
\begin{lemma}
The category $\mathbf{C}$ is homotopy equivalent to $\mathcal{C}$, as well the category $\mathbf{D}$ is homotopy equivalent to $\mathcal{D}$.
\end{lemma}

\subsubsection{Our strategy to solve it}
We say that an extension in $\EXT^k(A,A)$ is {\it degenerate} if it is the Yoneda product of two extensions $\alpha\in \EXT^a(A,A)$ and $\beta\in\EXT^b(A,A)$ where both $a,b>0$, and (at least) one of the extensions $\alpha$ and $\beta$ represents the 0 class in the $\Ext$ group.
Denote the full subcategory of $\EXT^k(A,A)$ consisting of the degenerate extensions by $\EXT^k_{\ddeg}(A,A)$.

Define category $\mathbf{Ext}^k(A,A)$. Its objects are morphisms $X\to Y$ of extensions from $\EXT^k(A,A)$, and morphisms are the commutative diagrams
\begin{equation}
\xymatrix{
X\ar[r]\ar[d]&X^\prime\ar[d]\\
Y\ar[r]&Y^\prime}
\end{equation}
We can replace, from the homotopical point of view, the categories $\EXT^k(A,A)$ by the categories $\mathbf{Ext}^k(A,A)$, as follows from the following simple lemma.

\begin{lemma}
The category $\mathbf{Ext}^k(A,A)$ is homotopy equivalent, for any $k$, to the category $\EXT^k(A,A)$.
\end{lemma}
\proof{}
\endproof

Define now the category $\mathbf{Ext}^k_{\ddeg}(A,A)$ as the full subcategory in $\mathbf{Ext}^k(A,A)$ whose objects are maps $X\to D$ where $D\in \EXT^k(A,A)$ is {\it degenerate}. Analogously to the lemma above, we prove
\begin{lemma}
The category $\mathbf{Ext}^k_\ddeg(A,A)$ is homotopy equivalent, for any $k$, to the category $\EXT^k_\ddeg(A,A)$.
\end{lemma}

Let $\mathcal{C}$ be an $n$-fold monoidal abelian category satisfying the condition (*). Then the category $\bigsqcup_k\EXT^k(A,A)$ is an $(n+1)$-fold monoidal category by ????. Moreover, it induces naturally an $(n+1)$-fold monoidal structure on the category $\bigsqcup_k\mathbf{Ext}^k(A,A)$.

Now the idea is to work with spectra of pairs, replacing the wrong-behaved ``base point'' by the subcategory $\mathbf{Ext}^k_{\ddeg}(A,A)$. It is motivated by the following lemma.

\begin{lemma}
Let $\mathcal{C}$ be an $n$-fold monoidal abelian category satisfying the condition (*).
Then any of $n+1$ monoidal structures on $\bigsqcup_k\mathbf{Ext}^k(A,A)$ defines a map of pairs
$(\mathbf{Ext}^a(A,A),\mathbf{Ext}^a_{\ddeg}(A,A))\times (\mathbf{Ext}^b(A,A),\mathbf{Ext}^b_{\ddeg}(A,A))\to (\mathbf{Ext}^{a+b}(A,A),\mathbf{Ext}^{a+b}_{\ddeg}(A,A))$.
The latter means that the image of $\mathbf{Ext}^a_{\ddeg}(A,A)\times \mathbf{Ext}^b(A,A)\sqcup \mathbf{Ext}^a(A,A)\times \mathbf{Ext}^b_{\ddeg}(A,A)$ for any of $n+1$ monoidal structure belongs to $\mathbf{Ext}^{a+b}_{\ddeg}(A,A)$.
\end{lemma}

\subsubsection{Chain homotopy of extensions}
S.Schwede [S] gives the following definition.
\begin{defn}
Let
$$
0\rightarrow N\rightarrow F_1\rightarrow \dots\rightarrow F_k\rightarrow M\rightarrow 0
$$
and
$$
0\rightarrow N\rightarrow G_1\rightarrow \dots\rightarrow G_k\rightarrow M\rightarrow 0
$$
be two extensions in $\EXT^k(M,N)$, and let $f,g\colon F^\mb\to G^\mb$ be two morphisms. We say that $f$ and $g$ are chain homotopic with fixed ends
if there is a sequence of maps $s_i\colon F_i\to G_{i-1}$ defined for $2\le i\le k$ such that
\begin{equation}
f_i-g_i=
\begin{cases}
s_2\circ d &\text{ for }i=1\\
d\circ s_i+s_{i+1}\circ d &\text{ for }2\le i\le k-1\\
d\circ s_k &\text{ for }i=k
\end{cases}
\end{equation}
\end{defn}

\begin{lemma}
Let $X,Y\in\EXT^k(M,N)$ be two extensions, and let $f,g\colon X\to Y$ be chain homotopic with fixed ends morphisms.
Then the loop
\begin{equation}
\xymatrix{
X\ar[r]^{f}&Y&X\ar[l]_{g}
}
\end{equation}
is homotopic to 0 (in the sense of the classifying spaces).
\end{lemma}

See [S], Lemma 4.2 for a proof.

\subsubsection{The most degenerate based point}

\endcomment

\subsubsection{}
The idea is to replace each category $\EXT^k(A,A)$ by {\it a pair} of a category and its subcategory, which has the same homotopy groups, as well as possess the $(n+1)$-fold monoidal structure on the corresponding disjoint unit. That is, we replace the based objects (which do not exist in a way compatible with monoidal structures) by based subcategories.

For any category $\EXT^k(A,A)$ define the category $\Omega_\free\EXT^k(A,A)$ as follows. Its objects are ``free loops'' of length 1
\begin{equation}
X\rightarrow Y \leftarrow X
\end{equation}
and the morphisms are the natural commutative diagrams.

This category has a full subcategory $\overline{\Omega}_\free\EXT^k(A,A)$ consisting of
\begin{equation}\label{ident}
X\xrightarrow{id}X\xleftarrow{id}X
\end{equation}
clearly this category is isomorphic to $\EXT^k(A,A)$.

We claim the following:
\begin{klemma}\label{pari}
All homotopy groups of the pair $\pi_i(\Omega_\free\EXT^k(A,A),\overline{\Omega}_\free\EXT^k(A,A))$ are isomorphic to $\pi_{i+1}(\EXT^k(A,A))$.
\end{klemma}

We give a proof of this Lemma a bit later. Now let us explain how it can help us.

Define for that another category of free double loops, $\Omega^2_\free\EXT^k(A,A)$. Its object is a diagram
\begin{equation}\label{doubleloop}
\xymatrix{
X\ar[r]\ar[d]&Y\ar[d]&X\ar[l]\ar[d]\\
Z\ar[r]&W&Z\ar[l]\\
X\ar[r]\ar[u]&Y\ar[u]&X\ar[l]\ar[u]}
\end{equation}
and the morphisms are natural commutative diagrams.

This category contains a full subcategory $\overline{\Omega}^2_\free\EXT^k(A,A)$ forming from diagrams \eqref{doubleloop} such that either the vertical or the horizontal lines are of the form \eqref{ident}.

Fix now a monoidal structure on $\bigsqcup\EXT^k(A,A)$. Then it defines a products on pairs of categories
\begin{equation}
\begin{aligned}
\ &(\Omega_\free\EXT^a(A,A),\overline{\Omega}_\free\EXT^a(A,A))\times (\Omega_\free\EXT^b(A,A),\overline{\Omega}_\free\EXT^b(A,A))\to\\
&\to(\Omega_\free^2\EXT^{a+b}(A,A),\overline{\Omega}_\free^2\EXT^{a+b}(A,A))
\end{aligned}
\end{equation}
Thus, at this point we achieved our goal for introduction of ``based points'' in the setting, which is necessary for the using of smash-products of spectra.

We can analogously define the categories $\Omega_\free^d\EXT^k(A,A)$ and $\overline{\Omega}_\free^d\EXT^k(A,A)$.
We have

\begin{klemma}\label{pariingeneral}
The homotopy groups of the pair $\pi_i(\Omega_\free^d\EXT^k(A,A),\overline{\Omega}_\free^d\EXT^k(A,A))$ is isomorphic for any $d$
to $\pi_{i+d}(\EXT^k(A,A))$.
\end{klemma}

\subsection{A proof of Key-Lemma \ref{pari}}
\subsubsection{A topological counterpart}
Before proving the Key-Lemma \ref{pari} let us explain what does it mean topologically.
Let $X$ be a topological space, and let $\Omega_\free X$ be its free loop space. It has a complicated homotopical type.
Nevertheless, one has
\begin{equation}
\pi_i(\Omega_\free X)=\pi_i X \oplus \pi_i \Omega X
\end{equation}
where $\Omega X$ is the based loop space.

Indeed, there is a fibration $p\colon \Omega_\free X\to X$ with the fiber $\Omega X$. It has a canonical section, sending $x\in X$ to the constant loop based at $x$. Therefore, the long exact sequence of the fibration splits, which gives the result.

Now we prove
\begin{prop}
$\pi_i(\Omega_\free X, X)=\pi_i \Omega X$ for any $i$.
\end{prop}
\proof{}
Consider the embedding $i\colon X\to \Omega_\free X$.
According to general principles of the Eckmann-Hilton duality,
\begin{equation}\label{EH}
\pi_i(\Omega_\free X, X)=\pi_{i-1}P
\end{equation}
where $P$ is the homotopy fiber of the embedding $i$. Let us compute this fiber.
Its point is a path in $\Omega_\free X$ from two fixed points belonging to the image $i(X)$. That is, it is a based 2-sphere in $X$, or an element of $\Omega^2X$.
We get from \eqref{EH}:
\begin{equation}
\pi_i(\Omega_\free X,X)=\pi_{i-1}\Omega^2X=\pi_i\Omega X=\pi_{i+1} X
\end{equation}
\endproof

\subsubsection{A categorical counterpart}

\subsection{Passage to spectra}
\comment
We will be very sketch here, remaining the proofs for the next paper.

Define a topological spectrum $\mathfrak{X}_k$ and its subspectrum $\mathfrak{X}_k^0$, as follows.

There are canonical maps $\Omega_\free^n\EXT^{k+n}(A,A)\to\Omega_\free^{n+1}\EXT^{k+n+1}(A,A)$.
Define
\begin{equation}
\mathfrak{X}_k=\underrightarrow{lim}_nB(\Omega_\free^n\EXT^{k+n}(A,A))
\end{equation}
and
\begin{equation}
\mathfrak{X}_k^0=\underrightarrow{lim}_nB(\overline{\Omega}_\free^n\EXT^{k+n}(A,A))
\end{equation}
(here $B(-)$ is the classifying space functor).

Now let $\mathcal{C}$ be an abelian $n$-fold monoidal category satisfying the condition (*).

Then the $n+1$ monoidal structures are well-defined on the pair of spectra $(\mathfrak{X},\mathfrak{X}^0)$ {\it up to homotopy}.
Moreover, these structures are compatible with the structure maps in spectra up to homotopy.

So the idea is like that: before we had everything strict, but we had not the based points which made impossible to use the smash-products of spectra. Now we replaced the spectrum $\{\EXT^k(A,A)\}$ by the pair of spectra $(\{\mathfrak{X}_k\},\{\mathfrak{X}_k^0\})$ which is ``based'', but we pay for that a big price; namely all maps are now defined up to homotopy.

Consider for any topological space $Y$ the suspension spectrum $\Sigma^\infty Y$. Then we have the operad of spectra $\{\Sigma^\infty D_k^n\}$ where $\{D_k^n\}$ is the $n$-dimensional little discs operad. Then the operad $[\Sigma^\infty Y,\Sigma^\infty D_k^n]_\mb$ is an operad of graded vector spaces.
Our goal is to prove the following theorem:

\begin{theorem}
The operad of graded vector spaces $[\Sigma^\infty Y,\Sigma^\infty D_k^n]_\mb$ acts, for any topological space $Y$, on the graded vector space
$[\Sigma^\infty Y, (\mathfrak{X},\mathfrak{X}^0)]_\mb$.
\end{theorem}
A proof is left for the sequel paper.

Now consider the case when $Y$ is the disjoint union of two points, that it, when $\Sigma^\infty Y$ is the sphere spectrum $\mathbb{S}$.
Then it is well-known that $[\mathbb{S},(\mathfrak{X},\mathfrak{X}^0)]_{-i}=\pi_i(\mathfrak{X},\mathfrak{X}^0)$.

Now we prove
\begin{lemma}
For any $i$
\begin{equation}
\pi_i(\mathfrak{X},\mathfrak{X}^0)=\pi_i(\EXT(A,A))=\Ext^i(A,A)
\end{equation}
\end{lemma}
\proof{}
\endproof

On the other hand
\begin{equation}
[\Sigma^\infty Y,\Sigma^\infty D_k^n]_{-i}=\underrightarrow{lim}[\Sigma^{\ell+i}Y,\Sigma^\ell D_k^n]
\end{equation}
In particular, in the case when $Y$ is the disjoint union of two points, we get $\pi_i^{stab}(D_k^n)$.
\endcomment
Consider the following ``spectrum of categories'' $\mathfrak{X}$:
$$
\mathfrak{X}_0=\EXT^0(A,A),\ \ \mathfrak{X}_1=\Omega_\free\EXT^2(A,A),\ \ \dots,\ \ \mathfrak{X}_n=\Omega_\free^n\EXT^{2n}(A,A),\ \ \dots
$$
First of all, we should explain in which sense $\mathfrak{X}$ is a spectrum.

One has the Retakh's map $G\colon \EXT^k(A,A)\to\Omega\EXT^{k+1}$, see Section 1.3.5. Applying this map twice, and the embedding $\Omega\hookrightarrow \Omega_\free$, we get a map
$\sigma\colon\mathfrak{X}_k\to\Omega_\free\mathfrak{X}_{k+1}$. These are the structure maps of our spectrum $\mathfrak{X}$.

The disjoint union of categories $\bigsqcup \mathfrak{X}_k$ is an $(n+1)$-monoidal category in a natural way, when $\mathcal{C}$ is an abelian $n$-fold monoidal category satisfying the condition (*).

Moreover, there is a subspectrum $\mathfrak{X}^0\subset\mathfrak{X}$, where
$$
\mathfrak{X}^0_n=\overline{\Omega}_\free^n\EXT^{2n}(A,A)
$$
such that each monoidal structure is a map
\begin{equation}
\otimes_i\colon (\mathfrak{X}_k,\mathfrak{X}^0_k)\times (\mathfrak{X}_\ell,\mathfrak{X}^0_\ell)\to(\mathfrak{X}_{k+\ell},\mathfrak{X}^0_{k+\ell})
\end{equation}

Now is the question: how these monoidal structures are compatible with the spectrum structure maps $\sigma$?
More precisely, consider the diagram

\begin{equation}\label{diagramx}
\xymatrix{
(\mathfrak{X}_k,\mathfrak{X}_k^0)\times(\mathfrak{X}_\ell,\mathfrak{X}_\ell^0)\ar[r]^{\otimes_i}\ar[d]_{\sigma\times id}&(\mathfrak{X}_{k+\ell},\mathfrak{X}_{k+\ell}^0)\ar[d]^{\sigma}\\
(\Omega_\free(\mathfrak{X}_k),\Omega_\free(\mathfrak{X}_k^0))\times(\mathfrak{X}_\ell,\mathfrak{X}_\ell^0)\ar[r]^{\otimes_i}&(\Omega_\free(\mathfrak{X}_{k+\ell}),
\Omega_\free(\mathfrak{X}_{k+\ell}^0))}
\end{equation}

\begin{lemma}
The diagram \eqref{diagramx}, and the analogous diagram for $\sigma$ applied to the right factor on the left arrow, is for any $i$ homotopically commutative. The latter means that after applying of the classifying space functor, the corresponding diagram of topological spaces is commutative up to homotopy.
\end{lemma}
The proof will be given in a later version of the paper.
\endproof

Denote by $B(-)$ the classifying space functor, and by $\Sigma^\infty Y$ the suspension spectrum of a topological space $Y$. Recall that we denote by $\{D_k^n\}_k$ the $n$-dimensional little disc operad.

Based on the Lemma above, we prove the following theorem:
\begin{theorem}
Let $\mathcal{C}$ be an $n$-fold monoidal abelian category satisfying (*).
Then there is a map in the homotopical category of spectra of topological spaces
\begin{equation}\label{themainmap}
\Sigma^\infty D_k^{n+1}\wedge \underbrace{(B(\mathfrak{X}),B(\mathfrak{X}^0))\wedge \dots \wedge (B(\mathfrak{X}),B(\mathfrak{X}^0))}_{k\text{ factors}}\to (B(\mathfrak{X},B(\mathfrak{X}^0)))
\end{equation}
which obeys the operad action equations.
\end{theorem}
Here $-\wedge-$ is the smash-product on the homotopical category of spectra, see [A], Part III, Lecture 4.

\proof{}
The statement follows from the Lemma above and from the Theorem from [BFSV] that the classifying space operad of the operad
$\mathcal{M}_{n}(d)$ (see Section 3.5) is homotopically equivalent to the n-dimensional little disc operad.
\endproof

\subsection{Applications}
All applications go the following line. Let $Y$ be a topological space; consider the suspension spectrum $\Sigma^\infty Y$.
This spectrum is a ``coalgebra object'' in the homotopical category. Here we mean that there are coassociative maps
\begin{equation}\label{coalgebra}
\Sigma^\infty Y\to \underbrace{\Sigma^\infty Y \wedge \Sigma^\infty Y\wedge\dots\wedge \Sigma^\infty Y}_{k\text{ factors}}
\end{equation}
for any $k\ge 1$.

The maps \eqref{coalgebra}, together with \eqref{themainmap}, give maps
\begin{equation}\label{y}
[\Sigma^\infty Y,\Sigma^\infty D_k^{n+1}]\otimes_\mathbb{Z}[\Sigma^\infty Y,(B(\mathfrak{X}),B(\mathfrak{X}^0))]^{\otimes k}\to
[\Sigma^\infty Y,(B(\mathfrak{X}),B(\mathfrak{X}^0))]
\end{equation}
Here all maps of spectra are graded abelian groups, and $\{[\Sigma^\infty Y,\Sigma^\infty D_k^{n+1}]\}_{k\ge 1}$ is an operad of graded abelian groups.

We have proved

\begin{theorem}
Let $Y$ be a based topological space. Then the operad $\{[\Sigma^\infty Y,\Sigma^\infty D_k^{n+1}]_\mb\}_{k\ge 1}$ of graded abelian groups acts on the graded abelian group $[\Sigma^\infty Y,(B(\mathfrak{X}),B(\mathfrak{X}^0))]_\mb$.
\end{theorem}

When $Y$ is the disjoint union of two points, the maps from the spectrum $\Sigma^\infty Y$ are just stable homotopy groups.
Thus, we get

\begin{theorem}
Let $\mathcal{C}$ be an abelian $n$-fold monoidal category satisfying the condition (*).
Then the operad $\{\oplus_i\pi_i^{stab}(D^{n+1}_k)\}_{k\ge 1}$ of abelian groups acts on the graded space $\oplus_j\Ext^j_{\mathcal{C}}(A,A)$.
\end{theorem}
\proof{}
One only needs to compute the stable homotopy groups $\pi_i^{stab}(B(\mathfrak{X}),B(\mathfrak{X}^0))$ of the pair. But the latter is equal to
$\Ext^i_\mathcal{C}(A,A)$ by Key-Lemma \ref{pari} and by the Retakh's theory, see Section 1.
\endproof

The Hurewicz homomorphism gives an isomorphism
$$
\pi_i^{stab}(X)\otimes_\mathbb{Z}\mathbb{Q}\eqto H_i(X,\mathbb{Q})
$$
where we get the link with the more usual statement that the homology operad of the little disc operad acts on the Hochschild
cohomology of any associative algebra.

In particular, for any bialgebra $A$ we constructed a 2-fold monoidal abelian category of tetramodules.
We know from Section 4.3 that the category of tetramodules obeys the condition (*) when $A$ is a Hopf algebra.
By the result of R.Taillefer (see Section 6) we get:

\begin{corr}
The Gerstenhaber-Schack cohomology of any Hopf algebra defined over $\mathbb{Q}$ is naturally a 3-algebra.
\end{corr}

\section{Appendix: The Gerstenhaber-Schack cohomology, after R.Taillefer}
Here we give an of overview of the works of R.Taillefer on Gerstenhaber-Schack cohomology, from slightly different point of view.
The main result is that the Gerstenhaber-Schack cohomology $H^\mb_\GS(A)$ (see Subsection 6.1) is equally equal to
$\Ext^\mb_{\Tetra(A)}(A,A)$ (Theorem 6.1 below). Let us notice that this result is true for any bialgebra, not necessarily a Hopf algebra.

\subsection{The Gerstenhaber-Schack complex}
Let $A$ be a (co)associative bialgebra.  Note that the
bar-differential in $\BBar^{\boxtimes_1}(A)$ is given by maps of
tetramodules; analogously, the cobar-differential in
$\Cobar^{\boxtimes_2}(A)$ is given by maps of tetramodules.

Let us recall, that originally the Gerstenhaber-Schack complex was
defined in [GS] as

\begin{equation}\label{gs}
C_{\GS}^\mb(A)=\Hom_{\Tetra(A)}(\BBar^{\boxtimes_1}_-(A),\Cobar^{\boxtimes_2}_+(A))
\end{equation}

Here $\BBar_-(B)$ and $\Cobar_+(C)$ are truncated complexes, which
end (start) with $B{\boxtimes_1}B$ ($C{\boxtimes_2}C$)
correspondingly.

For convenience of the reader let us write down here the
Gerstenhaber-Schack differential in $C_{\GS}^\mb(A)$ explicitly:

First of all, as a graded vector space,
\begin{equation}\label{eqeqeq_1}
C_{\GS}^\mb(A)=\oplus_{m,n\ge 0}\Hom_k(A^{\otimes m}, A^{\otimes
n})[-m-n]
\end{equation}
Now let $\Psi\colon A^{\otimes m}\to A^{\otimes n}\in
C_{GS}^{m+n}(A)$. We are going to define the Gerstenhaber-Schack
differential $d_\GS (\Psi)\in \Hom (A^{\otimes (m+1)},A^{\otimes
n})\oplus\Hom (A^{\otimes m},A^{\otimes (n+1)})$. Denote the
projection of $d_\GS$ to the first summand by $(d_\GS)_1$, and the
projection to the second summand by $(d_\GS)_2$. The formulas for
$(d_\GS)_1$ and $(d_\GS)_2$ are:
\begin{equation}\label{eq0.6fin}
\begin{aligned}
\ &(d_\GS)_1(\Psi)(a_0\otimes\dots\otimes a_m)=\ \ \ \ \ \ \ \ \ \ \ \ \ \ \ \ \ \ \ \ \ \ \hspace{7cm}\\
&\Delta^{n-1}(a_0)*\Psi(a_1\otimes\dots\otimes a_m)\\
&+\sum_{i=0}^{m-1}(-1)^{i+1}\Psi(a_0\otimes\dots\otimes (a_i*
a_{i+1})\otimes
\dots\otimes a_m)\\
&+(-1)^{m-1}\Psi(a_0\otimes\dots\otimes a_{m-1})*\Delta^{n-1}(a_m)
\end{aligned}
\end{equation}
and
\begin{equation}\label{eq0.7fin}
\begin{aligned}
\ &(d_\GS)_2(\Psi)(a_1\otimes\dots\otimes a_m)=\\
&(\Delta^{(1)}(a_1)*\Delta^{(1)}(a_2)*\dots
*\Delta^{(1)}(a_m))\otimes
\Psi(\Delta^{(2)}(a_1)\otimes\dots\otimes \Delta^{(2)}(a_m))\\
&+\sum_{i=1}^n(-1)^i\Delta_{i}\Psi(a_1\otimes\dots\otimes a_m)\\
&+(-1)^{n+1}\Psi(\Delta^{(1)}(a_1)\otimes\Delta^{(1)}(a_2)\otimes\dots
\otimes\Delta^{(1)}(a_m))\otimes(\Delta^{(2)}(a_1)*\Delta^{(2)}(a_2)*\dots
*\Delta^{(2)}(a_m))
\end{aligned}
\end{equation}

The goal of this Appendix is to prove the following Theorem due to R.Taillefer:

\begin{theorem}([Ta1,2])
Suppose a bialgebra $A$ has unit and counit. Then one has:
$$
\Ext^\mb_{\Tetra(A)}(A,A)=H^\mb\bigl(\Hom_{\Tetra(A)}(\BBar^{\boxtimes_1}_{-}(A),\Cobar^{\boxtimes_2}_{+}(A))\bigr)
$$
\end{theorem}

\subsection{Two forgetful functors and their adjoint}
\subsubsection{}
Let $A$ be a (co)associative bialgebra. Besides the category
$\Tetra(A)$, we can consider the categories $\Bimod(A)$ of
$A$-bimodules (when we consider $A$ as an algebra) and $\Bicomod(A)$
of $A$-bicomodules (when we consider $A$ as a coalgebra). Clearly
there are two exact forgetful functors $F_1\colon \Tetra(A)\to
\Bicomod(A)$ and $F_2\colon \Tetra(A)\to \Bimod(A)$. We have the
following
\begin{lemma}
Let $A$ be a bialgebra which has unit and counit. Then the functor
$F_1$ admits a left adjoint $L$ and the functor $F_2$ admits a right
adjoint $R$. The functors $L$ and $R$ are exact.
\end{lemma}
\begin{defn}
Let $A$ be a bialgebra, and let $N$ be an $A$-bicomodule, and let
$M$ be an $A$-bimodule.  The tetramodule $L(N)$ is called the
induced (from $N$) tetramodule, and the tetramodule $R(M)$ is called
the coinduced (from $M$) tetramodule. The induced and coinduced
tetramodules form full additive subcategories in the abelian
category $\Tetra(A)$. We denote them $\Tetra_\Ind(A)$ and
$\Tetra_\Coind(A)$, respectively.
\end{defn}

{\it Proof of Lemma:} we set
\begin{equation}\label{adjoint1}
L(N)=A\boxtimes_1 N\boxtimes_1 A
\end{equation}
and
\begin{equation}\label{adjoint2}
R(M)=A\boxtimes_2 M\boxtimes_2 A
\end{equation}
(see Section 4.2.1 for the definitions of $\boxtimes_1$ and
$\boxtimes_2$). Strictly speaking, to write formulas like this, $M$
and $N$ should be tetramodules. But the reader probably have noticed
that in the definition of $M_1\boxtimes_1 M_2$ we do not use the
right $A$-module structure in $M_1$ and the left $A$-module
structure in $M_2$. Analogously, in the definition of
$M_1\boxtimes_2 M_2$ we do not use the right comodule structure in
$M_1$ and the left comodule structure in $M_2$. Therefore,
(\ref{adjoint1}) and (\ref{adjoint2}) make sense.

We should check the adjunction properties
\begin{equation}\label{adjoint3}
\Hom_{\Bicomod(A)}(N,F_1(T))=\Hom_{\Tetra(A)}(L(N),T)
\end{equation}
and
\begin{equation}\label{adjoint4}
\Hom_{\Bimod(A)}(F_2(T),M)=\Hom_{\Tetra(A)}(T,R(M))
\end{equation}
as bifunctors.

Prove (\ref{adjoint3}). Any map of tetramodules
$L(N)=A\boxtimes_1N\boxtimes_1A\to T$ is uniquely defined by its
restriction to $1\boxtimes_1 N\boxtimes_11$. This map clearly is a
map of bicomodules $N\to F_1(T)$. Wise versa, any map of bicomodules
$N\to F_1(T)$ can be uniquely extended to a map of bimodules
$A\boxtimes_1 N\boxtimes_1 A\to T$ which is in fact a map of
tetramodules. These two assignments are inverse to each other.

The proof of (\ref{adjoint4}) is analogous.

The exactness of $L$ and $R$ is clear from their constructions
(\ref{adjoint1}) and (\ref{adjoint2}).
\endproof

\subsubsection{}
\begin{lemma}
Let $A$ be an associative bialgebra. Then any tetramodule
$M\in\Tetra(A)$ can be imbedded onto a coinduced tetramodule, and
there is a surjection into $M$ from an induced tetramodule.
\end{lemma}
\proof{} Let $M$ be an $A$-tetramodule. Consider
$P(M)=A\boxtimes_1 M\boxtimes_1 A$, it is induced from the
bicomodule $F_1(M)$. The map $p\colon A\boxtimes_1 M\boxtimes_1 A\to
M$, $a\boxtimes_1m\boxtimes_1b\mapsto a\cdot m\cdot b$ is clearly a
map (and an epimorphism, because $A$ contains a unit) of
tetramodules. Analogously, the tetramodule $Q(M)=A\boxtimes_2
M\boxtimes_2 A$ is coinduced from the bimodule $F_2(M)$, and we have
a monomorphism $j\colon M\to A\boxtimes_2 M\boxtimes_2 A$, $m\mapsto
\Delta_\ell\circ\Delta_r(m)$.
\endproof

\begin{corr}([Ta2])
For any bialgebra $A$ the category $\Tetra(A)$ has enough
injectives. \end{corr}
\proof{} The functor $R$ is a right adjoint to
an exact functor $F_2$, and, therefore, maps injective objects to
injective (see [W], Prop. 2.3.10). Moreover, it is left exact ([W],
Section 2.6). Therefore, it is sufficient to imbed $M$ as a
$A$-bimodule onto an injective $A$-bimodule $I$ (see [W], Section
2.3) and apply the functor $R$ to this this imbedding of
$A$-bimodules. This will give an imbedding $j\colon M\to Q(M)$ where
$Q(M)$ is defined in the proof of Lemma above.
\endproof

The main fact about the induced and the coinduced tetramodules is
the following

\begin{prop}
Let $A$ be a bialgebra with unit and counit. Then the functor
$X\mapsto \Hom_{\Tetra(A)}(X,Q)$ for fixed $Q\in\Tetra_{\Coind}(A)$
is an exact functor from $\Tetra_\Ind(A)^\opp$ to $\Ab$. As well,
the functor $Y\mapsto\Hom_{\Tetra(A)}(P,Y)$ for fixed
$P\in\Tetra_\Ind(A)$ is an exact functor from $\Tetra_\Coind(A)$ to
$\Ab$. \end{prop}
\proof{} Let us prove the first statement. Let
\begin{equation}\label{adjoint5}
0\rightarrow LN^\prime\rightarrow LN\rightarrow
LN^{\prime\prime}\rightarrow 0
\end{equation}
be an exact sequence of tetramodules,
$N,N^{\prime},N^{\prime\prime}\in \Bicomod(A)$. We should prove that
the sequence
\begin{equation}\label{adjoint6}
0\rightarrow\Hom_{\Tetra(A)}(LN^{\prime\prime},RM)\rightarrow\Hom_{\Tetra(A)}(LN,RM)\rightarrow\Hom_{\Tetra(A)}(LN^\prime,RM)\rightarrow
0
\end{equation}
is exact for any $M\in \Bimod(A)$.

By the adjunction, the exactness of (\ref{adjoint6}) is equivalent
to the exactness of the sequence
\begin{equation}\label{adjoint7}
0\rightarrow\Hom_{\Bimod(A)}(F_2LN^{\prime\prime},M)\rightarrow\Hom_{\Bimod(A)}(F_2LN,M)\rightarrow\Hom_{\Bimod(A)}(F_2LN^\prime,M)\rightarrow
0
\end{equation}
But this sequence is exact because for any $N\in\Bicomod(A)$ the
$A$-bimodule $F_2LN$ is free and, therefore, projective.

The second statement is proven analogously.
\endproof

\subsection{Some homological algebra}
In this Subsection we consider some homological algebra, which is
useful for computation of $Ext$ functors.
\subsubsection{A $(\mathcal{P},\mathcal{Q})$-pair}
Suppose $\mathcal{A}$ is an abelian category, and $\mathcal{P}$,
$\mathcal{Q}$ are two additive subcategories. We say that they form
a $(\mathcal{P},\mathcal{Q})$-pair if the following conditions are
satisfied:
\begin{itemize}
\item[1.] the functor $\Hom(?,Q)$ is exact on $\mathcal{P}^\opp$ for
any $Q\in\mathcal{Q}$;
\item[2.] the functor $\Hom(P,?)$ is exact on $\mathcal{Q}$ for any
$P\in\mathcal{P}$;
\item[3.] for any object $M\in\mathcal{A}$ there is an epimorphism
$P\to M$ for $P\in\mathcal{P}$ and there is a monomorphism $M\to Q$
for $Q\in\mathcal{Q}$;
\item[4.] a stronger version of 3: for any short exact sequence
$0\to M_1\to M_2\to M_3\to 0$ in $\mathcal{A}$ the epimorphisms
$P_i\to M_i$, $p_i\colon P_i\in\mathcal{P}$, can be chosen such that
there is a map of complexes
\begin{equation}\label{delta}
\xymatrix{
0\ar[r]&P_1\ar[r]\ar[d]^{p_i}&P_2\ar[r]\ar[d]^{p_2}&P_3\ar[r]\ar[d]^{p_3}&0\\
0\ar[r]&M_1\ar[r]&M_2\ar[r]&M_3\ar[r]&0}
\end{equation}
where the upper line is an exact sequence; also the dual condition
for $\mathcal{Q}$ for monomorphisms $q_i\colon M_i\to Q_i$,
$Q_i\in\mathcal{Q}$.
\end{itemize}
Note that the third condition guarantees that each object
$M\in\mathcal{A}$ has a $\mathbb{Z}_{\le 0}$-graded resolution in
$\mathcal{P}$ and a $\mathbb{Z}_{\ge 0}$-graded resolution in
$\mathcal{Q}$.
\begin{example}
If $\mathcal{A}$ has enough projectives and $\mathcal{P}$ is the
additive subcategory of projective objects,
$\mathcal{Q}=\mathcal{A}$ gives a $(\mathcal{P},\mathcal{Q})$-pair.
Analogously, if $\mathcal{A}$ has enough injectives and
$\mathcal{Q}$ is the additive subcategory of injective objects,
$\mathcal{P}=\mathcal{A}$ gives a $(\mathcal{P},\mathcal{Q})$-pair.
\end{example}
\begin{prop}
Suppose $\mathcal{A}=\Tetra(A)$ for an associative unital (see
Notations) bialgebra $A$. Then the pair
$(\Tetra_\Ind(A),\Tetra_\Coind(A))$ is a
$(\mathcal{P},\mathcal{Q})$-pair.
\end{prop}
\proof{} The first two properties were proven in Proposition 6.7, and
the third property follows from Lemma 6.5. Moreover, this
construction in this Lemma gives the fourth assertion in the
definition of a $(\mathcal{P},\mathcal{Q})$-pair immediately.
\endproof

\subsubsection{The Key-Lemma}
The main fact about $(\mathcal{P},\mathcal{Q})$-pairs is the
following lemma:
\begin{klemma}
Let $\mathcal{A}$ be an abelian category, having enough projective
or injective objects. Suppose we are given a
$(\mathcal{P},\mathcal{Q})$-pair and let $M,N\in\mathcal{A}$ be two
objects. Suppose $P^\mb\to M$ be a resolution of $M$ by objects in
$\mathcal{P}$, and $N\to Q^\mb$ be a resolution of $N$ by objects in
$\mathcal{Q}$. Then
\begin{equation}\label{kl}
\Ext_{\mathcal{A}}^\mb(M,N)=H^\mb(\Hom_{\mathcal{A}}(P^\mb,Q^\mb))
\end{equation}
\end{klemma}
\proof{} The proof consists from several steps. We give the proof
for the case of enough injectives, the case of enough projectives is
analogous. We recall the universal property which a derived functor
obeys, and prove that the functor $(M,N)\mapsto
H^\mb(\Hom_{\mathcal{A}}(P^\mb,Q^\mb))$ has this universal property.

Let $\mathcal{A},\mathcal{B}$ be  abelian categories, and $F\colon
\mathcal{A}\to\mathcal{B}$ be a left exact functor. Then the right
derived functor $\mathrm{L}^\mb F$ enjoys the following universal
property. We say that a collection of functors
$\{T_n\colon\mathcal{A}\to\mathcal{B}\}$, $n\ge 0$, is a
(cohomological) $\delta$-functor if for any exact sequence
\begin{equation}\label{ses}
0\to M\to N\to L\to 0
\end{equation}
in $\mathcal{A}$ one has a morphism $\delta\colon T_n(L)\to
T_{n+1}(M)$, $n\ge 0$, in $\mathcal{B}$ with the following long
exact sequence
\begin{equation}\label{les}
\dots\rightarrow T_{n-1}(L)\xrightarrow{\delta}T_n(M)\rightarrow
T_n(N)\rightarrow T_n(L)\xrightarrow{\delta}
T_{n+1}(M)\rightarrow\dots
\end{equation}
$n\ge 1$, depending functorially on the short exact sequence
(\ref{ses}). Consider a $\delta$-functor $\{T_n\}$. We say that this
$\delta$-functor is {\it universal} if for any other
$\delta$-functor $\{S_n\}$ with the natural transformation
$f_0\colon T_0\to S_0$ there is a unique morphism of
$\delta$-functors $\{f_n\colon T_n\to S_n\}$ extending $f_0$. From
this definition it follows that the universal $\delta$-functor with
$T_0=F$, if it exists, is unique. This point of view, independent on
existence of enough projective objects, due to Grothendieck [Tohoku]
and is extracted by the author from [W], Chapter 2.

Now we consider the $\Hom_{\mathcal{A}}(M,?)$ as a functor of the
second argument. If $\mathcal{A}$ has enough injectives the functors
$T_n(M,N)=\Ext^n_{\mathcal{A}}(M,N)$ form a cohomological universal
$\delta$-functor (see [W], Theorem 2.4.7).

Now our proof of the Key-Lemma will go as follows:
\begin{itemize}
\item[Step 1.] $T_k\colon(M,N)\mapsto H^k(\Hom_{\mathcal{A}}(P^\mb,Q^\mb))$
with $P^\mb\in\mathcal{P}$, $Q^\mb\in\mathcal{Q}$ is well-defined,
that is does not depend on the choice of $P^\mb$ and $Q^\mb$;
\item[Step 2.] it is a homological $\delta$-functor with
$T_0=\Hom_{\mathcal{A}}(M,N)$;
\item[Step 3.] it is a universal $\delta$-functor.
\end{itemize}
Clearly the Key-Lemma follows from these 3 statements. Let us prove
them.

{\it Step 1:} it simple follows from the conditions 1) and 2) in the
definition of a $(\mathcal{P},\mathcal{Q})$-pair.

{\it Step 2:} it follows easily from condition 4) in the definition
of a $(\mathcal{P},\mathcal{Q})$-pair.

{\it Step 3:} this is a bit more tricky. An additive functor
$F\colon \mathcal{A}\to\mathcal{B}$ is called {\it effaceable} if
for any object $N\in\mathcal{A}$ there is a monomorphism $j\colon
M\to I$ such that $F(j)=0$. It is proven in [Tohoku] that a
cohomological $\delta$-functor $\{T_n\}$ for which all $T_n$ for
$n\ge 1$ are effaceable, is universal. It remains to prove that our
functors $T_n(N)=H^n(\Hom_{\mathcal{A}}(P^\mb(M),Q^\mb(N)))$, $n\ge
1$, are effaceable. We can choose a monomorphism $j\colon N\to I$
with $I\in\mathcal{Q}$ by the condition 3) in the definition of a
$(\mathcal{P},\mathcal{Q})$-pair. Now the effaceability follows from
$T_n(I)=0$, $n\ge 1$, which immediately follows from 1) and 2).

Thus, it is proven that the functors
$H^n(\Hom_{\mathcal{A}}(P^\mb,Q^\mb))$ are universal
$\delta$-functors with the same 0-component than
$\Ext^n_{\mathcal{A}}(M,N)$. Therefore, these two functors coincide
because the universal functor $\{T_n\}$ with fixed $T_0$ is unique.

\endproof

\subsection{Example: computation of the Gerstenhaber-Schack
cohomology for $A=S(V)$} Here we compute, as an application of the
previous results of this Section, the Gerstenhaber-Schack cohomology
for the (co)free commutative and cocommutative bialgebra $A=S(V)$,
where $V$ is a finite-dimensional vector space. We prove
\begin{prop}
Let $A=S(V)$ be (co)free commutative cocommutative bialgebra, $V$
finite-dimensional. Then the Gerstenhaber-Schack cohomology is
$H^k_\GS(A)=\bigoplus_{i+j=k, i,j\ge
0}\Lambda^iV\otimes\Lambda^j(V^*)$. \end{prop}

\proof{}
We compute $\Ext^\mb_{\Tetra(A)}(A,A)$ for $A=S(V)$ as $H^\mb(\Hom(P^\mb, Q^\mb))$, where $P^\mb$ is a resolution of $A$ by free modules, and $Q^\mb$ is a resolution of $A$ by cofree comodules. Then we take the usual Koszul resolutions of the diagonal for $P^\mb$ and $Q^\mb$.
\endproof

As was noticed in Example 4.4, there is a canonical 3-algebra structure on the graded space $\oplus_{i+j=k}\Lambda^i V\otimes \Lambda^j V^*$.
It would be interesting to check explicitly that this 3-algebra structure coincides with the one defined in Section 4.4 via the 2-fold monoidal structure on the category of tetramodules.

\comment
\section{Appendix: explicit definition of 2-monoidal bialgebra}
In this short Appendix we write down explicitly what is a 2-monoidal bialgebra defined in Example \ref{monalg}.
The author believes that this is a correct 2-categorical generalization of the concept of bialgebra, in particular, that many features of the Drinfeld's theory of quantum groups exist in this context.

We just give a definition, and hope to consider this concept in more detail in the future.

{\it A 2-monoidal bialgebra} is a vector space $A$ with an associative algebra structure, denoted by $*\colon A^{\otimes 2}\to A$, and with two coassociative coproducts, denoted by $\Delta_1,\Delta_2\colon A\to A^{\otimes 2}$, which are compatible.

So far we have the following axioms:
\begin{itemize}
\item[1)] $(a*b)*c=a*(b*c)$;
\item[2)] $(\Delta_1\otimes 1)\circ\Delta_1 (a)=(1\otimes \Delta_1)\circ \Delta_1(a)$;
\item[3)] $(\Delta_2\otimes 1)\circ\Delta_2 (a)=(1\otimes \Delta_2)\circ \Delta_2(a)$;
\item[4)] $\Delta_1(a*b)=\Delta_1(a)*\Delta_1(b)$;
\item[5)] $\Delta_2(a*b)=\Delta_2(a)*\Delta_2(b)$.
\end{itemize}

To formulate the remaining axioms, recall that we use the notation $\Delta_i(a)=\Delta_i^1(a)\otimes \Delta_i^2(a)$ which assumes
$\Delta_i(a)=\sum_j\Delta_i^{j1}(a)\otimes`\Delta_i^{j2}(a)$, where $i=1$ or $i=2$.

\begin{itemize}
\item[6)] for any $a\in A$, there is a linear map $\eta(a)\colon (\Delta_2\circ \Delta_1^1(a))\otimes (\Delta_2\circ \Delta_1^2(a))\to
(\Delta_1\circ \Delta_2^1(a))\otimes (\Delta_1\circ \Delta_2^2(a))$ which satisfy the following axioms 7) and 8) below;
\item[7)] for any $a\in A$ one has:
\begin{equation}
\xymatrix{
(\Delta_2\Delta_1^1\Delta_1^1(a))\otimes (\Delta_2\Delta_1^2\Delta_1^1(a))\otimes (\Delta_2\Delta_1^2(a))\ar[r]^{\eta(\Delta_1^1(a))\otimes id_{A\otimes A}}&(\Delta_1\Delta_2^1\Delta_1^1(a))\otimes(\Delta_1\Delta_2^2\Delta_1^1(a))\otimes (\Delta_2\Delta_1^2(a))\ar[d]^{\eta(a)}\\
ghglihilhiljh&(\Delta_1\otimes \Delta_1\otimes id_{A\otimes A})\bigl((\Delta_1\Delta_2^1(a))\otimes (\Delta_1\Delta_2^2(a))\bigr)}
\end{equation}
\item[8)]
\end{itemize}

In particular, the roles of $\Delta_1$ and $\Delta_2$ are not symmetric.

\endcomment

\bigskip

\noindent {\sc Faculty of Science, Technology and Communication,
Campus Limpertsberg, University of Luxembourg, 162A avenue de la
Faiencerie, Luxembourg, \\L-1511 LUXEMBOURG}

\bigskip

\noindent {\em E-mail address\/}: {\tt borya$\_$port@yahoo.com}

\end{document}